\documentclass[11pt, letterpaper]{article}
\usepackage[utf8]{inputenc}
\usepackage{geometry}
\usepackage[spanish, english, es-nodecimaldot]{babel}

\usepackage{graphicx}
\usepackage{latexsym, mathtools, amssymb, amsfonts, amsbsy}
\usepackage[backend=bibtex8, citestyle=authoryear, backref, natbib, style=authoryear, sorting=nyt]{biblatex}
\DefineBibliographyStrings{english}{%
  backrefpage = {\lowercase{s}ee p.},%
  backrefpages = {\lowercase{s}ee pp.}%
}
\usepackage[backref=pages]{hyperref}
\hypersetup{colorlinks,
			citecolor=DeepPink2,
			filecolor=black,
			urlcolor=black,
			linkcolor=DeepPink2}
\usepackage{enumitem}
\usepackage[x11names]{xcolor}
\usepackage[english, noabbrev]{cleveref}
\usepackage{relsize}
\usepackage{bbm}
\usepackage{float}
\usepackage{caption}
\usepackage{units}
\usepackage{subfigure}

\usepackage{algorithm}
\usepackage{algpseudocode}
\algrenewcommand{\algorithmiccomment}[1]{\textcolor{blue}{/*#1*/}}

\title{On the Consistency of a Random Forest Algorithm in the Presence of Missing Entries}
\author{Irving Gómez-Méndez \and Emilien Joly}
\date{\today}

\newtheorem{hypothesis}{Hypothesis}
\newtheorem{theorem}{Theorem}
\newtheorem{proposition}{Proposition}
\newtheorem{lemma}{Lemma}
\newtheorem{techlemma}{Technical Lemma}

\newcommand{\indicadora}[2]{\mathlarger{\mathbbm{1}}\,_{#1}#2}

\let\tinymatrix\smallmatrix

\patchcmd{\tinymatrix}{\vcenter}{\vtop}{}{}
\patchcmd{\tinymatrix}{\bgroup}{\bgroup\scriptsize}{}{}

\usepackage{tikz}

\usetikzlibrary{babel}
\usetikzlibrary{matrix, automata}
\usetikzlibrary{arrows, arrows.meta}
\usetikzlibrary{shadows}
\usetikzlibrary{positioning}
\usetikzlibrary{trees}
\usetikzlibrary{decorations.pathmorphing}
\definecolor{olivegreen}{rgb}{0,0.6,0}
\usetikzlibrary{calc}
\usetikzlibrary{shapes.multipart}
\usetikzlibrary{matrix, fit}
\usetikzlibrary{backgrounds}
\tikzset{point/.style={insert path={ node[scale=3*sqrt(\pgflinewidth)]{.} }}}
\tikzset{pointgreen/.style={insert path={ node[scale=3*sqrt(\pgflinewidth), color = SpringGreen4]{.} }}}
\tikzset{pointblue/.style={insert path={ node[scale=3*sqrt(\pgflinewidth), color = blue]{.} }}}
\tikzset{pointred/.style={insert path={ node[scale=3*sqrt(\pgflinewidth), color = red]{.} }}}
\tikzset{Import_line/.style n args={3}{very thick,color=#1!#2!#3},
	Import_line/.default={black}{50}{black}}

\definecolor{purples1}{RGB}{252,251,253}
\definecolor{purples2}{RGB}{239,237,245}
\definecolor{purples3}{RGB}{218,218,235}
\definecolor{purples4}{RGB}{188,189,220}
\definecolor{purples5}{RGB}{158,154,200}
\definecolor{purples6}{RGB}{128,125,186}
\definecolor{purples7}{RGB}{106,81,163}
\definecolor{purples8}{RGB}{84,39,143}
\definecolor{purples9}{RGB}{63,0,125}

\def\X{\mathbf{X}}

\def\B{\mathbf{B}}
\def\N{\mathbb{N}}
\def\M{\mathbf{M}}
\def\x{\mathbf{x}}

\def\imiss{\mathbf{i}}

\def\P{\mathbb{P}}          % Probabilidad
\def\R{\mathbb{R}}          % Numeros reales
\def\V{\mathbb{V}}          % Varianza
\def\E{\mathbb{E}}          % Esperanza
\def\sample{\mathcal{D}_n}

\newcommand{\llavs}[1]{\left\lbrace #1 \right\rbrace}%       {}
\newcommand{\corchet}[1]{\left\lbrack #1 \right\rbrack}%     []
\newcommand{\parent}[1]{\left( #1 \right)}%                  ()
\newcommand{\abs}[1]{\lvert#1\rvert}%                    ||
\newcommand{\norm}[1]{\lVert#1\rVert}%                 || ||

\DeclareMathOperator*{\argmin}{arg\,min}
\DeclareMathOperator*{\argmax}{arg\,max}

\begin{document}

\maketitle

\begin{abstract}
This paper tackles the problem of constructing a non-parametric predictor when the latent variables are given with incomplete information. The convenient predictor for this task is the random forest algorithm in conjunction to the so-called CART criterion. The proposed technique enables a partial imputation of the missing values in the data set in a way that suits both a consistent estimator of the regression function as well as a partial recovery of the missing values.  A proof of the consistency of the random forest estimator is given in the case where each latent variable is missing completely at random (MCAR).
\end{abstract}

\section{Introduction}

Random forests algorithms are widely used in applied computer science nowadays. They provide very useful tools for extracting information of a possibly large scaled data set throughout two tasks: regression and classification. First introduced by \citet{breiman2001random}, this non-parametric technique has various computational benefits over many regression/prediction techniques as it is fast and simple. Therefore it is really common to see the random forest algorithm (or parts of it) as a building block of more intricate machine learning algorithms. The small number of hyper-parameters to be tuned, along with the bootstrap-aggregation technique \citep{breiman1996bagging} -commonly called bagging- allows a straightforward parallelization and has made random forests one of the most popular tools for handling data sets of large size.
It has been successfully involved in various practical problems, including chemioinformatics \citep{svetnik2003random}, ecology \citep{prasad2006newer,cutler2007random}, 3D object recognition \citep{shotton2011real}, bioinformatics \citep{diaz2006gene} and econometrics \citep{varian2014big}.

On the theoretical side, little is known about the mathematical guarantees of the method. This has lead the community to underline a ``gap'' between theory and practice. Part of this gap can be explained by the bagging mechanism and the splitting criterion. Each of these processes introduces a source of randomness into the construction of the trees which makes the random forest algorithms very challenging to study in their full generality.
One way to overcome the theoretical difficulties of random forests is through simplified versions of the original procedure. This is often done by simply ignoring the bagging step and/or by replacing the splitting criterion with a more elementary splitting protocol (e.g. \citet{breiman2004consistency}). However, in recent years, important theoretical studies have been performed to analyze more elaborated models (e.g. \citet{biau2008consistency,ishwaran2010consistency,biau2012analysis,genuer2012variance,zhu2015reinforcement}). Consistency and asymptotic normality for Breiman's infinite forests were proven in \citet{wager2018estimation} simplifying the splitting step. \citet{mentch2016quantifying} proved a similar result for finite forests and \citet{scornet2015consistency} proved a consistency result in the context of additive regression models.
The review presented by \citet{biau2016random} gives a very good introduction to those concepts and proposes solutions for some of the challenges presented above.

The relative simplicity of random forests advocates for its use to handle missing data. As an initial solution, \citet{Breiman2003}
proposed an algorithm which takes advantage of the so-called proximity-matrix to impute the missing values and \citet{ishioka2013imputation} modified it in order to obtain a new version more robust to outliers. The algorithm MissForest, proposed by \citet{stekhoven2011missforest} is another alternative to iteratively improve the imputation of the missing values using random forests. An interesting alternative comes from algorithms that attempt to handle the missing values  without performing an imputation. The so-called surrogate splits (see \citet{breiman1984CART,venables2002,ripley2007pattern,friedman2009elements,hapfelmeier2012recursive}) is an example of this family of algorithms. This technique has spread from the original CART trees and has been introduced in conditional trees \citep{hothorn2006unbiased}. Missing Incorporated in Attributes (MIA) \citep{twala2008good} is another approach that has driven attention in recent years. In a practicle point of view, \citet{josse2019consistency} perform a simulation study to compare different methods to handle missing data using trees and study the consistency of two approaches to estimate the prediction function with missing values in a framework of multiple imputation \citep{rubin2004multiple,rubin1996multiple}. An overview of different methodologies to handle missing data can be found in \citet{josse2018introduction}.
%, the use of a universally consistent algorithm and for a Missing At Random (MAR) mechanism.% Many of these algorithms has been implemented in \texttt{R} in the package \texttt{partykit} \citep{JMLR:v16:hothorn15a} 

In the present paper, we propose a new splitting protocol and prove the consistency of the algorithm in the context of an additive model in presence of missing data accordingly to a Missing Completely at Random (MCAR) paradigm. Moreover, in a simulation study, we compare this new algorithm with other approaches that handle missing values using random forests also built with the CART criterion. We choose the following regression function  which has been used in previous simulation studies \citep{breiman1996bagging,rieger2010random,josse2019consistency,friedberg2020local}, \[m(\x)=10\sin\parent{\pi\x^{(1)}\x^{(2)}}+20\parent{\x^{(3)}-0.5}^2+10\x^{(4)}+5\x^{(5)}.\]
Following the schema presented by \citet{rieger2010random}, we simulate i.i.d. observations uniformly distributed on $[0,1]^5$ and introduce missing values in $\X^{(1)},\X^{(3)}$ and $\X^{(4)}$. We let the percentage of the missing values in $\X^{(4)}$ vary, to show the performance of this new algorithm in a possibly highly corrupted sample.

The paper is organized as follows. In \Cref{sec:RF_aggregating}, we introduce the random forests, the concept of data-missing mechanisms and show the problems presented by the original CART criterion when there are missing entries. The main consistency result is presented in \Cref{sec:main_results} and its proof is given in \Cref{sec:Proofs}. In \Cref{sec:simulations}, we describe the simulation study and present the conclusions of this study.

\section{Random Forests: Aggregating Decision Trees}
\label{sec:RF_aggregating}
The random forests algorithms consist in building a collection of objects called decision trees (which are rooted trees) that, most of the time, are binary trees. For a specific tree, each path starting from the root corresponds to a decision (sometimes called prediction) and the descent mechanism is usually a key part in the decision. Each bifurcation in a single tree is done using a random component and an optimization component, representing a split of the space into two parts. Then, moving along a path in the decision tree represents a decision between one of the two parts of the space. Finally, the random forest is made of a set of decision trees that are later aggregated all together. This aggregation step is also important. In this paper, this takes the form of a simple mean (see \Cref{eq:random_forest} for a precise statement). The main strength of random forests is that it aggregates the information of many different decision trees in a global predictor that ends to be a lot more stable (and then a lot more informative) than each specific tree.

Besides being conceptually simple, decision trees are powerful and attractive in practice for several reasons:

\begin{itemize}
\item They can model arbitrarily complex relations between the input and the output space.
\item They handle categorical or numerical variables, or a mix of both.
\item They can be used in regression or supervised classification problems.
\item They are easily interpretable, even for non-statisticians.
\end{itemize}

\begin{figure}[ht]
\begin{center}
\begin{minipage}[t]{0.3\textwidth}
	\begin{tikzpicture}[scale=1.7]
	% Dibujamos el espacio X = [0,1]x[0,1]
	\draw[Import_line] (0,0) rectangle (1,1);

	\draw[Import_line, fill=purples6] (0,0.5) rectangle (1,1);
	\draw[Import_line, fill=purples5] (0,0) rectangle (1,0.5);

		% Y el nombre a los ejes
	\node (X1) at (0.5,0)[below] {$\mathsmaller{\X^{(1)}}$};
	\node (X2) at (0,0.5) [left] {$\mathsmaller{\X^{(2)}}$};

		% Los cortes
	\node (z1) at (1,0.5) [right] {$\color{red} z_1$};
	\draw[Import_line, color=red] (0,0.5) -- (1,0.5);

		% Aqui estan los puntos
	\draw (0.85,0.15) node[point] {};
	\draw (0.65,0.35) node[point] {};
	\draw (0.10,0.65) node[point] {};
	\draw (0.70,0.85) node[point] {};
	\draw (0.25,0.25) node[point] {};
	\end{tikzpicture}

	\begin{tikzpicture}[grow=down,shape=circle,very thick,level distance=18mm]
	\tikzstyle{level 1}=[sibling distance=20mm]

	\tikzstyle{edge from parent}=[color=black,draw]

	\node [label={right:\small $\X^{(2)}\geq \color{red} z_1$}, style = {fill = purples6}] {}
		child {node [style = {fill = purples6}] {}}
		child {node [style = {fill = purples7}]{}};
	\end{tikzpicture}
\end{minipage}
\hspace{0.1cm}
\begin{minipage}[t]{0.3\textwidth}
	\begin{tikzpicture}[scale=1.7]
	% Dibujamos el espacio X = [0,1]x[0,1]
	\draw[Import_line] (0,0) rectangle (1,1);
	\draw[Import_line, fill=purples6] (0,0.5) rectangle (1,1);
	\draw[Import_line, fill=purples6] (0,0) rectangle (0.75,0.5);
	\draw[Import_line, fill=purples4] (0.75,0) rectangle (1,0.5);

	% Y el nombre a los ejes
	\node (X1) at (0.5,0)[below] {$\mathsmaller{\X^{(1)}}$};
	\node (X2) at (0,0.5) [left] {$\mathsmaller{\X^{(2)}}$};

		% Los cortes
	\node (z1) at (1,0.5) [right] {$z_1$};
	\node (z2) at (0.75,1) [above] {$\color{red} z_2$};
	\draw[Import_line, color=red] (0.75,0) -- (0.75,0.5);

		% Aqui estan los puntos
	\draw (0.85,0.15) node[point] {};
	\draw (0.65,0.35) node[point] {};
	\draw (0.10,0.65) node[point] {};
	\draw (0.70,0.85) node[point] {};
	\draw (0.25,0.25) node[point] {};
	\end{tikzpicture}

	\begin{tikzpicture}[grow=down,shape=circle,very thick,level distance=18mm]
	\tikzstyle{level 1}=[sibling distance=20mm]
	\tikzstyle{level 2}=[sibling distance=12mm]

		\tikzstyle{edge from parent}=[ color=black,draw]

		\node [label={right:\small $\X^{(2)}\geq z_1$}, style = {fill = purples6}] {}
		child {node [label={[label distance=0.1cm]350:\footnotesize $\X^{(1)}\geq \color{red} z_2$},style = {fill = purples6}] {}
			child {node [style = {fill = purples7}]{}}
			child {node [style = {fill = purples5}]{}}
			}
		child {node [style = {fill = purples7}]{}
			};
	\end{tikzpicture}
\end{minipage}
	\hspace{0.1cm}
\begin{minipage}[t]{0.3\textwidth}
	\begin{tikzpicture}[scale=1.7]
	% Dibujamos el espacio X = [0,1]x[0,1]
	\draw[Import_line] (0,0) rectangle (1,1);
	\draw[Import_line, fill=purples4] (0,0.5) rectangle (0.4,1);
	\draw[Import_line, fill=purples6] (0,0) rectangle (0.75,0.5);
	\draw[Import_line, fill=purples8] (0.4,0.5) rectangle (1,1);
	\draw[Import_line, fill=purples4] (0.75,0) rectangle (1,0.5);

	% Y el nombre a los ejes
	\node (X1) at (0.5,0)[below] {$\mathsmaller{\X^{(1)}}$};
	\node (X2) at (0,0.5) [left] {$\mathsmaller{\X^{(2)}}$};

		% Los cortes
	\node (z1) at (1,0.5) [right] {$z_1$};
	\node (z2) at (0.75,1) [above] {$z_2$};
	\node (z3) at (0.4,1) [above] {$\color{red} z_3$};
	\draw[Import_line, color=red] (0.4,0.5) -- (0.4,1);

		% Aqui estan los puntos
	\draw (0.85,0.15) node[point] {};
	\draw (0.65,0.35) node[point] {};
	\draw (0.10,0.65) node[point] {};
	\draw (0.70,0.85) node[point] {};
	\draw (0.25,0.25) node[point] {};
	\end{tikzpicture}

		\begin{tikzpicture}[grow=down,shape=circle,very thick,level distance=18mm]
	\tikzstyle{level 1}=[sibling distance=20mm]
	\tikzstyle{level 2}=[sibling distance=12mm]

		\tikzstyle{edge from parent}=[ color=black,draw]

		\node [label={right:\small $\X^{(2)}\geq z_1$}, style = {fill = purples6}] {}
		child {node [label={[label distance=0.1cm]350:\footnotesize $\X^{(1)}\geq z_2$},style = {fill = purples6}] {}
			child {node [style = {fill = purples7}]{}}
			child {node [style = {fill = purples5}]{}}
			}
		child {node [label={[label distance=0.1cm]350:\footnotesize $\X^{(1)}\geq \color{red} z_3$}, style = {fill = purples7}]{}
			child {node [style = {fill = purples5}]{}}
			child {node [style = {fill = purples9}]{}}
			};
	\end{tikzpicture}
\end{minipage}
\caption{At each step of the tree construction a partition is performed over a cell (or equivalently its corresponding node) maximizing some split criterion.}
\label{fig:trees}
\end{center}
\end{figure}
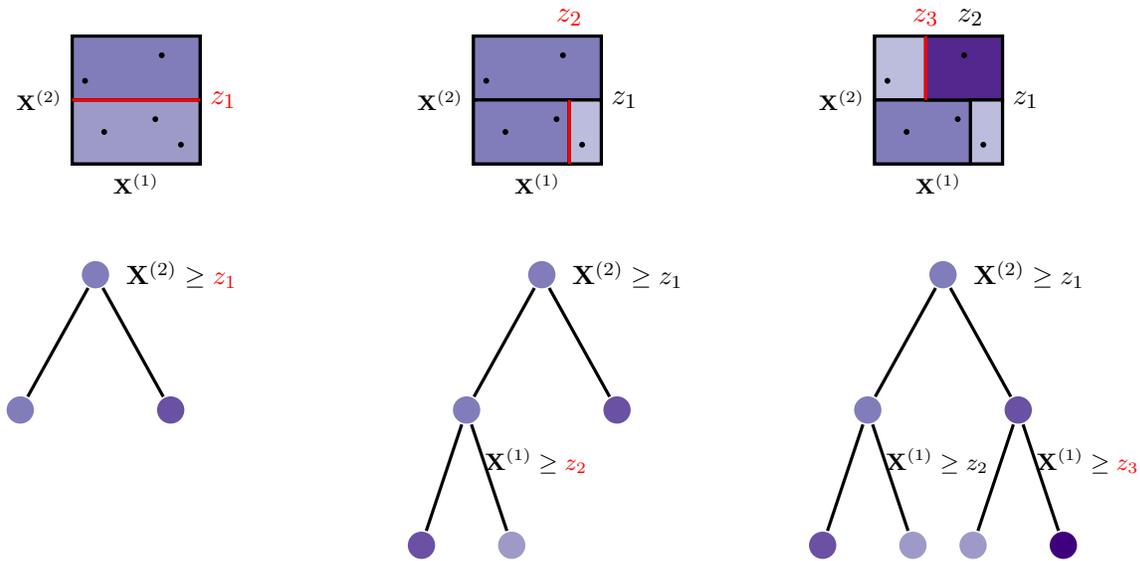

To construct the partition of the input space, decision trees work in a recursive way. The root of the tree is the whole input space, $\mathcal{X}$, which is split into disjoint regions. Then each region is split into more regions, and this process is continued until some stopping rule is applied (see \Cref{fig:trees}). Note that, for each cell considered during the construction of the tree, we can estimate the target variable. Hence, the decision tree is the collection of all the nodes and the estimation of the target variable for each one of those nodes. At each step of the tree construction, the partition performed over a cell (or equivalently its corresponding node) is determined by maximizing some splitting criterion (see Equation \eqref{eq:CART_criterion} for a proper definition).

In our framework, we assume to have access to a training set $\sample=(\X_i,Y_i)_{i=1,\dots,n}$ where the response variables $Y_i$ are real-valued and the input variables $\X_i$ belong to some space $\mathcal{X}$. For simplicity, we assume throughout the discussion that the vector \(\X\) is a uniform random variable on the square \(\mathcal{X}\). Although appearing restrictive at first sight, this assumption is not essential for the following. It could be replaced by an assumption of the form: $\X$ has a positive density over the full space $\mathcal{X}$. In most applications the space $\mathcal{X}$ is a compact portion of a $p$ dimensional space. Hence we assume, with very little loss of generality, that $\mathcal{X}=[0,1]^p$.
The task is to predict the random variable $Y$ with respect to the input vector $\X=(\X^{(1)},\ldots,\X^{(p)})$ where $\X^{(j)}\in [0,1]$ (for $j=1,\ldots,p$).
For example, one could be interested in finding a function $f:[0,1]^p\to\R$ that minimizes the functional loss $\E_{\X,Y}\corchet{\mathcal{L}(f(\X),Y)}$ where $\mathcal{L}:\R^2\to\R$ is the squared error $\mathcal{L}(f(\X),Y)=(f(\X)-Y)^2$. The solution of this optimization problem
\[m=\argmin_{f:[0,1]^p\to\R}\E_{\X,Y}\corchet{(f(\X)-Y)^2}\]
is given by the regression function $m(\X)=\E[Y|\X]$.
For practical problems, the distribution of $(\X,Y)$ and hence, the regression function, is unknown. The task is then to use the data $\sample$ to construct a learning model, also called learner, predictor or estimator, $m_n:[0,1]^p\to\R$ that estimates the regression function $m$, and enables us to predict the outcome for new incoming inputs.
At least two distinct strategies can be used to solve the estimation problem: the parametric and the non-parametric estimation. In a nutshell, the non-parametric estimation uses a set of functions \(\mathcal{F}_n\) (often called model) that serves the minimization (often referred to Empirical Risk Minimization)
\begin{equation*}
  m_n=\argmin_{f\in \mathcal{F}_n}\E_{\sample}\corchet{(f(\X)-Y)^2}.
\end{equation*}
The main focus of this work is to construct a non-parametric estimation of \(m\) through random forests algorithms.
As mentioned above, a random forest is a predictor consisting of $M(>1)$ randomized trees. The randomization is introduced in two different steps of the tree construction and will be governed by the following parameters. 
%Prior to the construction of each tree, $a_n$ observations are extracted at random with (or without) replacement from the learning data set $\sample$. Only these $a_n$ observations are taken into account in the tree construction. Then, at each cell (or node), a split is performed by maximizing the split criterion over a number $\texttt{mtry}$ of input variables chosen uniformly at random among the original ones. The tree construction is stopped when each final node contains less or equal than \texttt{nodesize} points or when the tree has $t_n$ final nodes. Hence, the parameters of this algorithm are:

\begin{itemize}

\item An integer $a_n\in\{1,\ldots,n\}$ defines the number of observations to be subsampled in $\mathcal{D}_n$ for the construction of each tree. Then, a tree built upon those observations (say $\X_{i_1},\dots,\X_{i_{a_n}}$) depends on $\mathcal{D}_n$ only through those variables.

\item An integer $\texttt{mtry}\in\{1,\ldots,p\}$ that corresponds to the number of directions (features) chosen randomly at each step. Those selected features are denoted by $\mathcal{M}_{try}\subset \{1,\ldots,p\}$ and are redefined at each node of the tree construction.

\item Two integers $\texttt{nodesize}$ and $q_n$ such that $2q_n-1\le \texttt{nodesize}\le a_n$. The number $\texttt{nodesize}$ (resp.  $q_n$) holds for the maximum (resp. minimum) number of observations in a final cell. In the algorithm, a cell $A$ that contains a number of observation $q_n\le N(A) \le \texttt{nodesize}$ is declared \textbf{final} and will not be split again. 

\end{itemize}

This randomization (independent from the original source of randomness in the sample \(\sample\)) is represented in a symbolic random variable \(\Theta\). To each tree -- randomized with the random variable \(\Theta_k\) -- is associated a predicted value at a query point $\x$, denoted as $m_n(\x;\Theta_k)$. The different trees are constructed by the same procedure but with independent randomization so the random variables $\Theta_1,\ldots,\Theta_M$ are i.i.d. with common law $\Theta$. The nature and dimension of $\Theta$ depends on its use in the tree construction. In our choice of construction rules, \(\Theta\) consists of the observations selected for the tree and the candidate directions to be split at each step.
The tree construction stops when all the remaining cells are final cells.
At the end of the tree construction, a partition of the space \(\mathcal{X}\) is returned in a form of the collection of final cells \((A_{n,i})_{i\ge 1}\). Each of these cells corresponds to a leaf of the tree.
Finally, the $k$-th tree estimate is defined as

\begin{equation}
m_n(\mathbf{x};\Theta_k)=\sum_{i\in \mathcal{I}_{n,\Theta_k}}\frac{Y_i\indicadora{\X_i\in A_n(\x;\Theta_k)}{}}{N(A_n(\x;\Theta_k))}, \label{eq:tree}    
\end{equation}

where $\mathcal{I}_{n,\Theta_k}$ is the set of the \(a_n\) observations selected prior to the construction of the $k$th tree, $A_n(\x;\Theta_k)$ is the unique final cell that contains $\x$, and $N(A_n(\x;\Theta_k))$ is the number of observations which belong to the cell $A_n(\x;\Theta_k)$.
The aggregation of the trees forms the finite random forest estimator given by

\begin{equation}
m_{M,n}(\x;\Theta_1,\ldots,\Theta_M)=\frac{1}{M}\sum_{k=1}^M m_n(\x;\Theta_k).\label{eq:random_forest}    
\end{equation}

It is known from the work of \citet{breiman2001random} that the random forest estimator does not overfit when $M$ tends to infinity so the value of \(M\) is only restricted by the computational power. Besides being important in practice, we do not make big case of the choice of \(M\) since the following results show the consistency of each of the tree estimators, hence the consistency of the random forest for any number \(M\).
It remains a challenge to be able to show the consistency of the random forest in the case when some trees are not consistent estimators.

%\[m_n(\x;\sample)=\E_{\Theta|\sample}[m_n(\x;\Theta,\sample)]\]

Different split criteria have been proposed depending on the statistical problem and the nature of the input space. For regression purposes, the most commonly used is the CART criterion \citep{breiman1984CART}. For supervised classification, the split criterion takes the form of an impurity function, like the misclassification error, the Gini index \citep{gini1912variabilita} or the Shannon entropy \citep{shannon1948mathematical}, considered in the criteria known as ID3 and C4.5 \citep{quinlan1986induction,quinlan1993}.

\paragraph{The CART criterion}
The present work focuses on the CART split criterion that we introduce now. We begin with some important notations.

\begin{itemize}
\item $A$ denotes a general node (or cell).

\item $N(A)$ holds for the number of points in $A$.

\item The notation $d=(h,z)$ denotes a cut in $A$, where
the integer $h\in\{1,\ldots,p\}$ is called a direction (or feature) and
$z \in [0,1]$ is the position of the cut in the $h$-th direction.

\item $\mathcal{C}_A$ is the set of all possible cuts in the node $A$. It means that \(h\) do belong to the set $\mathcal{M}_{try}$ and that for any chosen \(h\), \(z\) lies between the bounds of the cell \(A\) in that specific direction (see for example Figure \ref{fig:trees}).

\item A cell $A$ is split accordingly to a cut $d= (h,z)$ into two cells denoted as $A_L=\{\x\in A\,:\,\x^{(h)}< z\}$ and $A_R=\{\x\in A\,:\,\x^{(h)}\geq z\}$.

\item $\bar{Y}_A$ (resp. $\bar{Y}_{A_L}$, $\bar{Y}_{A_R}$) is the empirical mean of the response variable $Y_i$ for the indexes such that $\X_i$ belongs to the cell $A$ (resp. $A_L$, $A_R$).
\end{itemize}
Then, the empirical version of the CART split criterion for a generic cell $A$ is defined as
\begin{align}
L_n(A,d)=&\frac{1}{N(A)}\sum_{i=1}^n\parent{Y_i-\bar{Y}_A}^2\indicadora{\X_i\in A}{}\nonumber\\
&-\frac{1}{N(A)}\sum_{i=1}^n\left(Y_i-\bar{Y}_{A_L}\indicadora{\X_i^{(h)}<z}{}-\bar{Y}_{A_R}\indicadora{\X_i^{(h)}{}\geq z}\right)^2\indicadora{\X_i\in A}{},\label{eq:CART_criterion}
\end{align}
with the convention $0/0=0$.

Intuitively, the CART criterion compares (in a weighted manner) the empirical variance inside the cell \(A\) with the sum of the two empirical variances of the sub-cells \(A_R\) and \(A_L\). At each step of the creation of the random tree, a current cell \(A\) is selected and cut by choosing the best empirical cut $\widehat{d}=(\widehat{h},\widehat{z})$ so that $L_n(A,d)$ is maximal over the set $\mathcal{C}_A$, that is
\[
\widehat{d}=(\widehat{h},\widehat{z})\in\argmax_{\begin{tinymatrix}
d\in \mathcal{C}_A
\end{tinymatrix}} L_n(A,d).
\]

\paragraph{In the presence of missing values}
Let us emphasize, through an example, the issues that one faces with the original CART criterion in a context of missing values. Consider the space $\mathcal{X}=[0,1]^2$ and  two observations, $\X_1$ and $\X_2$,  that belong to a cell $A \subset \mathcal{X}$. We assume that a direction and location for a cut of $A$ have been chosen and we denote by $A_L$ and $A_R$ the two resulting cells (see \Cref{tab:miss_CART}, \Cref{fig:miss_CART} for an illustration).
Since its value is missing, $\X_1^{(1)}$ is represented as a dashed line %\footnote{In our illustrations we represent observations with missing values with dashed lines.}
over the interval $[0,1]$.

\begin{table}[ht]
\begin{minipage}[b]{0.45\textwidth}
\centering
\begin{tabular}{ccc}
       & $\x^{(1)}$ & $\x^{(2)}$ \\
\hline
$\X_1$   & \texttt{NA} & 0.5 \\
$\X_2$   & 0.75        & 0.25\\
\hline
$A$      & [0.3,0.9]  & [0.2,0.7]\\
$A_L$    & [0.3,0.6]  & [0.2,0.7]\\
$A_R$    & [0.6,0.9]  & [0.2,0.7]\\
\end{tabular}
\caption{Toy example data. We consider a cell $A$ that has been split into $A_L$ and $A_R$. The feature $\X^{(1)}$ is missing in the first observation (denoted as \texttt{NA}).}
\label{tab:miss_CART}
\end{minipage}
\begin{minipage}[t]{0.45\textwidth}
\begin{center}
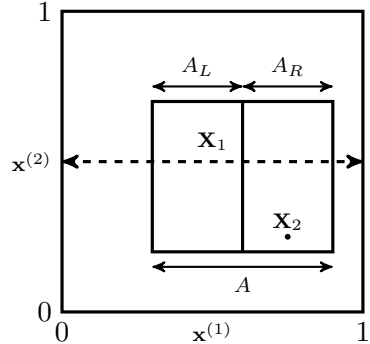

	\begin{tikzpicture}[scale=4, >=stealth']

	% Dibujamos el espacio X = [0,1]x[0,1]
	\draw[Import_line] (0,0) rectangle (1,1);

	% Y el nombre a los ejes
	\node (X1) at (0.5,0)[below] {$\mathsmaller{\x^{(1)}}$};
	\node (X2) at (0,0.5) [left] {$\mathsmaller{\x^{(2)}}$};

	% Limites a los ejes
	\draw (0,0) node[below] {0};
	\draw (0,0) node[left] {0};
	\draw (1,0) node[below] {1};
	\draw (0,1) node[left] {1};

	% Aqui estan los puntos
	\draw (0.75,0.25) node[point] {};
	\node at (0.75,0.3) {$\mathsmaller{\X_2}$};

	\draw[Import_line, dashed, <->] (0,0.5) -- node[above] {$\mathsmaller{\X_1}$} (1,0.5);

	% Aqui esta la celda
	\draw[Import_line] (0.3,0.2) rectangle (0.9,0.7);
	\draw[thick, <->] (0.3,0.15) -- node[below] {$\mathsmaller{A}$} (0.9,0.15);

	% Aqui esta el corte
	\draw[Import_line] (0.6,0.2) -- (0.6,0.7);

	% Las celdas hijas
	\draw[thick, <->] (0.3,0.75) -- node[above] {$\mathsmaller{A_L}$} (0.6,0.75);
	\draw[thick, <->] (0.6,0.75) -- node[above] {$\mathsmaller{A_R}$} (0.9,0.75);
	\end{tikzpicture}
\end{center}
\captionof{figure}{Illustration of the example data, $\X_1$ is represented as a dashed line since $\X^{(1)}$ is missing.}
\label{fig:miss_CART}
\end{minipage}
\end{table}

It is clear that $\X_2\in A_R$, however it is not possible to decide, without any further operation, if $\X_1\in A_R$ or $\X_1\in A_L$, or even if $\X_1\in A$.
The CART criterion is then intractable since the quantities $N(A)$, $N(A_L)$, $N(A_R)$, $\bar{Y}_A$, $\bar{Y}_{A_L}$, $\bar{Y}_{A_R}$, $\indicadora{\X_i\in A}{}$, $\indicadora{\X_i^{(h)}<z}{}$ and $\indicadora{\X_i^{(h)}\geq z}{}$ can not be computed. We propose a new partial imputation algorithm in Section \ref{sec:def_algo} to handle these problems.

\paragraph{Definition of Missingness}

The concept of data-missing mechanisms (introduced by \citet{rubin1976inference}) establishes the relationship between missingness and the data. It is common to define the data-missing mechanisms through the data matrix. However, to have a useful definition from a theoretical point of view, we formally define them using the random variables $\X$ and $Y$.  First,  let us define a new variable, called the indicator of missing value,
\[\M^{(h)}=\left\{\begin{array}{ll}
1 & \text{if } \X^{(h)}\text{ is missing}\\
0 & \text{otherwise}
\end{array} \right.,\quad 1\leq h\leq p.\]
A data point \((\X,Y)\) such that \(\M^{(h)}=1\) is said to be \textit{missing in the direction} \(h\). 
We are assuming throughout this work that the response $Y$ has no missing values, so it is not necessary to define an indicator of missing value for $Y$. The missingness mechanisms are fully characterized by the information of the conditional distribution of $\M^{(h)}$ given $(\X,Y)$.
One of the data-missing mechanisms identified by \citet{rubin1976inference} is the so-called Missing Completely at Random (MCAR). We say that the variable $\X^{(h)}$ is MCAR if $\M^{(h)}$ is independent of $(\X,Y)$. In other words, under the MCAR assumption, a coordinate $\X^{(h)}$ has some probability to be missing in the sample and this probability does not depend on the value of \(\X\) nor the response variable \(Y\). Note that it is allowed to get missing rates that differ from one direction to another (i.e. $(\M^{(h)})_h$ are not identically distributed) as seen in the hypothesis of Theorem \ref{T1}.

In this work, we allow the probability of missingness to depend on $n$. This hypothesis is relevant when the corrupted data in a sample represents a proportion of the entire sample; for example when measuring or sending data of fixed size in a network, a certain number of bits is lost, in average. Note that this does not imply that the distribution of the data is depending on the sample size $n$ but the fact that one specific entry is lost do depend on $n$. We define the expected proportion of missing values by $p_n^{(h)}=\P[\M^{(h)}=1|\sample]$. See \Cref{T1} for the condition on the values of \(p_n^{(h)}\).

\section{Main Result}
\label{sec:main_results}
\subsection{Introduction of our Algorithm}
\label{sec:def_algo}

The approach proposed in this paper overpasses the CART computation problems through iterative partial imputations of the missing values that are refined at each step. Unlike most of the so-called imputation techniques, the imputation step is \textit{not} performed independently of the evaluation of the CART criterion. On the contrary, the imputation is itself part of the optimization. The split of a cell (in the usual framework) is then replaced by the couple (split, imputation), where the additional imputation part refers to assign deterministically each data point missing in the direction of the split to one of the two resulting cells. We repeat this step on the resulting cells until a stopping rule is achieved. In the present article we consider splits that keep at least a determined number of points \(q_n\) to each child node, and we stop when there is no split with such property.

This algorithm looks similar to other approaches that perform imputation, except for the fact that the imputations take the form of closed intervals rather than punctual values. Each iteration of the algorithm consists of the (split, imputation) task throughout the optimization of the CART criterion. As a result, each iteration takes as an input the current cell \(A\) and the current imputations of the missing variables belonging to \(A\), and respond a pair of cells \((A_L,A_R)\) (for left and right) and new imputations. The initialization of the imputation of each missing variable is given by \(\widehat{\X}_i^{(h)}=[0,1]\) which is interpreted as giving to the missing value all the possible values in \([0,1]\). To avoid confusion in the sequel, we make the subtle difference between the imputation of the missing values at the start (referred to as \textit{in}) of an iteration and the imputation at the end (referred to \textit{out}) of the iteration. 

For one cut $d=(h,z)$, we give the following formalism for the assignations of the missing random variables in the direction $h$. An assignation is a binary vector $w$ of size equal to the number of missing values in the direction $h$ that belong to $A$, and such that $w_i=1$ if the $i$th missing observation is chosen to belong to the left cell $A_L$ and $w_i=0$ if it is chosen to belong to the right cell $A_R$. The result of this operation is a vector $\widehat{\X}_{w}$, such that
\[\widehat{\X}_{i,w}^{(h)}=\left\{\begin{array}{ll}
\corchet{a^{(h)},z} & \text{if } w_{i}=1\\
\corchet{z,b^{(h)}} & \text{if } w_{i}=0
\end{array} \right.\]
if the value of $\X_i$ is missing in the direction $h$ and $\widehat{\X}_{i,w}^{(h)}=\X_i^{(h)}$ otherwise. The values $a^{(h)}$ and $b^{(h)}$ hold for the extreme point of the cell $A$ in the direction $h$.
See \Cref{fig:example_proposal} for an illustration with $p=2$.
\begin{center}
\begin{figure}[ht]
\begin{minipage}[t]{0.45\textwidth}
\begin{center}
\begin{tikzpicture}[scale=4, >=stealth']
	% Dibujamos el espacio X = [0,1]x[0,1]
	\draw[Import_line] (0,0) rectangle (1,1);

	% El nombre de la celda
	\draw[thick, <->] (0,1.07) -- node[above] {$\mathsmaller{A}$} (1,1.07);

	% Y el nombre a los ejes
%	\node at (0.5,0)[below] {$\mathsmaller{\X^{(1)}}$};
%	\node at (0,0.5) [left] {$\mathsmaller{\X^{(2)}}$};

	% Aqui estan los puntos
	\node at (0.75,0.85)[above] {$\mathsmaller{\widehat{\X}_{1,in}}$};
	\draw (0.75,0.85) node[point] {};

	\node at (0.9,0.6)[above] {$\mathsmaller{\widehat{\X}_{3,in}}$};
	\draw (0.9,0.6) node[point] {};

	\node at (0.45,0.55)[above] {$\mathsmaller{\widehat{\X}_{4,in}}$};
	\draw (0.45,0.55) node[point] {};

	\node at (0.3,0.1)[above] {$\mathsmaller{\widehat{\X}_{7,in}}$};
	\draw (0.3,0.1) node[point] {};

	\node at (0.5,0.75)[above] {$\mathsmaller{\widehat{\X}_{2,in}}$};
	\draw[Import_line, dashed, <->] (0,0.75) -- (1,0.75);

	\node at (0.25,0.5)[above] {$\mathsmaller{\widehat{\X}_{5,in}}$};
	\draw[Import_line, dashed, <->] (0,0.5) -- (1,0.5);

	\node at (0.5,0.25)[above] {$\mathsmaller{\widehat{\X}_{6,in}}$};
	\draw[Import_line, dashed, <->] (0,0.25) -- (1,0.25);

	% Aqui esta el corte
	%\draw[Import_line, color=red] (0.6,0) -- (0.6,1);
\end{tikzpicture}
\end{center}
\end{minipage}
\hspace{0.1cm}
\begin{minipage}[t]{0.45\textwidth}
\begin{center}
\begin{tikzpicture}[scale=4, >=stealth']
	% Dibujamos el espacio X = [0,1]x[0,1]
	\draw[Import_line] (0,0) rectangle (1,1);

	% El nombre de la celda
	\draw[thick, <->] (0,1.07) -- node[above] {$\mathsmaller{A_L}$} (0.6,1.07);
	\draw[thick, <->] (0.6,1.07) -- node[above] {$\mathsmaller{A_R}$} (1,1.07);

	% Y el nombre a los ejes
%	\node at (0.5,0)[below] {$\mathsmaller{\X^{(1)}}$};
%	\node at (0,0.5) [left] {$\mathsmaller{\X^{(2)}}$};

	% Aqui estan los puntos
	\node at (0.75,0.85)[above] {$\mathsmaller{\widehat{\X}_{1,out}}$};
	\draw (0.75,0.85) node[point] {};

	\node at (0.9,0.6)[above] {$\mathsmaller{\widehat{\X}_{3,out}}$};
	\draw (0.9,0.6) node[point] {};

	\node at (0.45,0.55)[above] {$\mathsmaller{\widehat{\X}_{4,out}}$};
	\draw (0.45,0.55) node[point] {};

	\node at (0.3,0.1)[above] {$\mathsmaller{\widehat{\X}_{7,out}}$};
	\draw (0.3,0.1) node[point] {};

	\node at (0.3,0.75)[above] {$\mathsmaller{\widehat{\X}_{2,out}}$};
	\draw[Import_line, dashed, <->] (0,0.75) -- (0.6,0.75);

	\node at (0.8,0.5)[above] {$\mathsmaller{\widehat{\X}_{5,out}}$};
	\draw[Import_line, dashed, <->] (0.6,0.5) -- (1,0.5);

	\node at (0.3,0.25)[above] {$\mathsmaller{\widehat{\X}_{6,out}}$};
	\draw[Import_line, dashed, <->] (0,0.25) -- (0.6,0.25);

	% Aqui esta el corte
	\draw[Import_line, color=red] (0.6,0) -- (0.6,1);
\end{tikzpicture}
\end{center}
\end{minipage}
\caption{We perform a cut and assignation of points were the variable is missing, maximizing the CART criterion $L_n\parent{A,d,w}$.}\label{fig:example_proposal}
\end{figure}
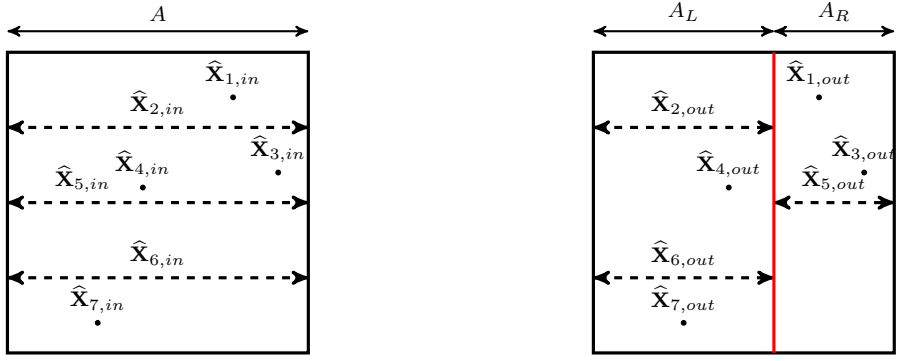
\end{center}
Finally, denoting by $I(A)$ the indexes of the observations belonging the cell $A$ (by previous assignation or being completely observed), the CART criterion in the context of missing values is defined on a cell $A$ by the formula
\begin{align*}
L_n\parent{A,d,w}=
&\frac{1}{\widehat{N}(A)}\sum_{i\in I(A)}\parent{Y_i-\widehat{Y}_{A}}^2\\
&-\frac{1}{\widehat{N}(A)}\sum_{i\in I(A)}\parent{Y_i-\widehat{Y}_{A_{L}}}^2\indicadora{a^{(h)}\le\widehat{\X}_{i,w}^{(h)}\le z}{}\\
&-\frac{1}{\widehat{N}(A)}\sum_{i\in I(A)}\parent{Y_i-\widehat{Y}_{A_{R}}}^2\indicadora{z<\widehat{\X}_{i,w}^{(h)} \le b^{(h)}}{}
\end{align*}
where $\widehat{Y}_{A_L}$ (resp. $\widehat{Y}_{A_R}$) is the empirical mean of the $Y_i$ such that $\widehat{\X}_{i,w}$ belongs to the cell $A_L$ (resp. $A_R$) and $\widehat{N}(A)$ is the number of points of the cell $A$.

Let $\widehat{N}_{miss}^{(h)}(A)=\sum_{i=1}^n \indicadora{\M_i^{(h)}=1}{}$ be the number of observations assigned to the cell $A$ that are missing in the direction $h$, note that $\sum_i w_i$ is the number of these observations which are assigned to the left node. Thus, $\widehat{p}_L^{(h)}=\sum_i w_i/\widehat{N}_{miss}^{(h)}(A)$ estimates the probability of an observation with the direction $h$ missing to belong to the left node, given that it belongs to the cell $A$. Similarly $\widehat{p}_R^{(h)}=1-\widehat{p}_L^{(h)}$ denotes the probability of an observation with the direction $h$ missing to belong to the right node, given that it belongs to the cell $A$. We can use these probabilities to predict the target variable for a new observation whose direction $h$ is missing. These can be done stochastically, assigning the new observation to the left node with probability $\widehat{p}_L^{(h)}$ and to the right node with probability $\widehat{p}_R^{(h)}$. If the new observation has a direction missing that was always observed during the training phase, then we might not be able to assign this observation to the left or right nodes when that particular direction was chosen in the construction of the tree. In such case, we can still estimate the target value using the estimation of the current cell. The prediction of the random forest would still be the average of the predictions of the trees. Note that this average would stabilize the stochastic behavior of the individual predictors, in a ``law of large numbers'' fashion way.

We give the precise algorithm for the CART optimization with assignation in Algorithm \ref{algo:CART_opti} and the construction of the random forest in Algorithm \ref{alg:random_forests_with_missingness}.

\begin{algorithm}[ht]
%\small
\hspace*{\algorithmicindent} \textbf{Input:} A cell $A$, the corresponding data $(\widehat{\X}_{i,in},Y_i)_i$ that belongs to $A$\\
\hspace*{\algorithmicindent} \textbf{Output:} Two cells $A_L,A_R$ and assignations $(\widehat{\X}_{i,out},Y_i)_i$ of the data.
\begin{algorithmic}[1]
\For{ each cut $d=(h,z)$ in $\mathcal{C}_A$}
	\State Compute $\widehat{N}_{miss}^{(h)}(A)=\sum_{i=1}^n \indicadora{\M_i^{(h)}=1}{}$.
	\State Let $\mathcal{W}_A^{(h)}=\{0,1\}^{\widehat{N}^{(h)}_{miss}(A)}$ be the possible assignations in the direction $h$.
	\State Let $\overline{\mathcal{W}_A^{(h)}}\subset \mathcal{W}_A^{(h)}$ such that $A_L$ and $A_R$ have at least $q_n$ elements.
	\If{$\overline{\mathcal{W}_A^{(h)}}=\emptyset$\label{alg:CART_check_empty}}
	    \State Remove $d$ from $\mathcal{C}_A$.
	\Else 
	    \State Compute $\widehat{w}_d\in\argmax_{
w\in\overline{\mathcal{W}_A^{(h)}}} L_n\parent{A,d,w}$.
	\EndIf
\EndFor
\State Compute $\widehat{d}\in\argmax_{
d\in\mathcal{C}_A} L_n\parent{A,d,\widehat{w}_d}$
\State Define $\widehat{\X}_{i,out}^{(h)}=\widehat{\X}_{i,\widehat{w}_{\widehat{d}}}^{(h)}$
\State Define $A_L =\{\x \in A:\x^{(\widehat{h})}\le \widehat{z}\}$ and $A_R =\{\x \in A:\x^{(\widehat{h})}> \widehat{z}\}$

\end{algorithmic}
\caption{CART optimization with assignation.}
\label{algo:CART_opti}
\end{algorithm}

One could ask if the condition for a possible (cut, assignation) in line \ref{alg:CART_check_empty} in Algorithm \ref{algo:CART_opti} could never be fulfilled. This case happens if and only if the number of data points in $A$ is less than $2q_n$. This condition is check in Algorithm \ref{alg:random_forests_with_missingness} anytime before Algorithm \ref{algo:CART_opti} is called. It basically checks if the resulting two cells in a split have enough number of points. Such condition is important to ensure the convergence of our algorithm and is restated in Theorem \ref{T1}. The task of Algorithm \ref{algo:CART_opti} is to compute the optimal couple (split, assignation) given by
\[(\widehat{d},\widehat{w})\in\argmax_{\begin{tinymatrix}
d\in \mathcal{C}_{A}\\
w\in\mathcal{W}_A^{(h)}
\end{tinymatrix}} L_n\parent{A,d,w}.\]

\begin{algorithm}[ht]
%\small
\hspace*{\algorithmicindent} \textbf{Input:} Training sample $\sample$, number of trees $M>1$, \texttt{mtry}$\in\{1,\ldots,p\}$, $a_n\in\{1,\ldots,n\}$, $\texttt{nodesize}\in\{2 q_n-1,\ldots,a_n\}$, $q_n\in\{1,\ldots,\lfloor a_n/2\rfloor\}$.\\
\hspace*{\algorithmicindent} \textbf{Output:} Random forest $m_{M,n}$.
\begin{algorithmic}[1]
\For{$i=1,\ldots,M$}
	\State Select $a_n$ points uniformly in $\sample$, denoted (abusively) $(\X_1,Y_1),\dots,(\X_{a_n},Y_{a_n})$.
	\State Set \(\widehat{\X}_i^{(h)}=[0,1]\) for each missing variable.
	\State Let $\mathcal{A}=\{[0,1]^p\}$, the initialization of the set of \textbf{active} cells.
	\State Let $m_n(\cdot,\Theta_i)=\emptyset$, the initialization of the decision tree.
	\While{$\mathcal{A}\neq \emptyset$}
	    \State Let $A$ be the first element of $\mathcal{A}$.
	    \If{$q_n \leq N(A) \leq \texttt{nodesize}$}
	        \State Remove $A$ from $\mathcal{A}$, add mark it as a \textbf{final} cell.
	        \State Compute the mean of the points in $A$, denoted as $\widehat{Y}_A$, and add $\{A, \widehat{Y}_A\}$ to $m_n(\cdot,\Theta_i)$.
	    \Else
    	    \State Select uniformly, without replacement, a subset $\mathcal{M}_{try}\subset \{1,\ldots,h\}$ of cardinality \texttt{mtry}.
    	    \State Apply \Cref{algo:CART_opti} over the directions in $\mathcal{M}_{try}$ and $\mathcal{C}_A$.
    	    \State Let be $\widehat{d}=(\widehat{h},\widehat{z})$ the optimal cut selected and $\widehat{w}$ the optimal assignation.
    	    \State Compute the mean of the points in $A$, denoted as $\widehat{Y}_A$.
    	    \State Compute the estimated probabilities $\widehat{p}_L^{(\widehat{h})}$ and $\widehat{p}_R^{(\widehat{h})}$
    	    \State Add $\left\{A, \widehat{Y}_A, \widehat{d},\widehat{p}_L^{(\widehat{h})},\widehat{p}_R^{(\widehat{h})}\right\}$ to $m_n(\cdot,\Theta_i)$.
    	    \State Remove $A$ from $\mathcal{A}$, and add $A_L$ and $A_R$ to $\mathcal{A}$.
        \EndIf
	\EndWhile
\EndFor 
\State Set $m_{M,n}=\{m_n(\cdot,\Theta_i)\}_i$.
\end{algorithmic}
\caption{Random forest with assignation of missing entries.} \label{alg:random_forests_with_missingness}
\end{algorithm}

\paragraph{On admissible assignations.}
At first sight, one could think that the number of possible assignations defined by $\mathcal{W}_A^{(h)}$ in Algorithm \ref{algo:CART_opti} contains all the binary vectors of size $\widehat{N}^{(h)}_{miss}(A)$ (that we denote abusively $N$ in this paragraph). This would lead to exponential complexity in the number of calculations of the CART criterion and then to an intractable algorithm. The following remark shows that only a linear (in $N$) number of assignations has to be considered.
Fix a cell $A$ and a cut $(h,z)$ in $A$, to keep a simple notation we denote by  $Y_{(1)}\le \ldots \le Y_{(N)}$ the ordered values of the response variable such that $\M^{(h)}=1$. Moreover, denote $\bar{Y}_{L,obs}$ (resp. $\bar{Y}_{R,obs}$) the empirical mean of the non-missing variables belonging to the cell $A_L$ (resp. $A_R$). Suppose, without lost of generality, that $\bar{Y}_{L,obs}\leq \bar{Y}_{R,obs}$. 
Since $\bar{Y}_{L,obs}\leq\bar{Y}_{R,obs}$, the only candidates to the maximization of the CART criterion are the one such that there exists a $0\le k \le N$ such that $Y_{(1)},\ldots,Y_{(k)} \in A_L$ and $Y_{(k+1)},\ldots,Y_{(N)} \in A_R$. The case $k=N$ (resp. $k=0$) holds for an assignation of all the variables to the left (resp. to the right). 
Indeed, if the assignation $w$ does not belong to those candidates, there exists a couple such that $Y_{(i)}< Y_{(j)}$ and $Y_{(i)} \in A_R$, $Y_{(j)} \in A_L$. The assignation obtained by the simple transposition of those two variables strictly increases the CART criterion so that those assignations cannot be optimal.
Finally, this shows that only $N+1$ assignation are admissible. The cases where no observed variables are in $A_L$ or in $A_R$ can be treated in the same way.

In particular, note that if we have extra information on the missingness mechanism, the admissible assignations $w$ might be even reduced. For example, assume that the only  missing values are in $\X^{(1)}$ and that those values are missing exactly when they are bigger than some threshold $\tau$, that is, when we have right censoring. Also consider that we know the value of $\tau$, so the missingness mechanism is completely known. So, if we perform a cut at $\tau$ in $\X^{(1)}$ then the only admissible assignation for $\widehat{\X}_i$ is to the right child node.

\subsection{Hypothesis and Main Theorem}

We consider an additive regression model satisfying the following properties,

\begin{hypothesis}[H\ref{H1}] \label{H1}
The response variable $Y$ is of the form
\[Y=\sum_{j=1}^pm_j(\X^{(j)})+\varepsilon,\]
where $\X$ is uniformly distributed over $[0,1]^p$, $\varepsilon$ is an independent Gaussian centered noise with finite variance $\sigma^2>0$ and each component $m_j$ is continuous.
\end{hypothesis}

\begin{hypothesis}[H\ref{H2}]\label{H2}
The random variables \(\X^{(h)}\) are missing by following an MCAR mechanism. The probability of missingness $p_n^{(h)}=\P\corchet{\M^{(h)}=1}$ only depends on the size \(n\) of the sample $\sample$ and $\lim_{n\to\infty}p_n^{(h)}=c^{(h)}$  where $0<c^{(h)}<1$ is constant for all $h\in\{1,\ldots,p\}$.
\end{hypothesis}

\begin{theorem}\label{T1}
Assume that H\ref{H1} and H\ref{H2} hold. Then, under the condition $q_n\to\infty$, the random forest estimator with missing values is consistent in probability, i.e., for all $\xi,\rho>0$ there exists $N\in\N^\star$, such that for all $n>N$
\[\P\Bigl[|m_n(\X)-m(\X)|\le\xi\Bigr]\ge 1-\rho\]
\end{theorem}

Additive regression models, which decompose the regression function as a sum of univariate functions, are flexible and easy to interpret, providing a good trade-off between model complexity, calculation time and interpretation. We have taken the previous work of \citet{scornet2015consistency} as inspiration, they proved the consistency of the original random forest algorithm (without missing values) considering the same structure for the regression function and the same assumptions for the input variables and the errors. However, the proof of the results differs from theirs in several parts. Considering other assumptions for the input variables and the errors, like  sub-Gaussian errors or distributions with support in $[0,1]^p$ for $\X$, would require to adapt the technical parts of the proof with little relevance.

The fact that $q_n\to\infty$ implies that the number of points selected in each tree $a_n$, tends to infinity. Hence, in the sequel, we assume that $a_n\to\infty$ and omit the dependence of the trees over $a_n$. The condition that $q_n\to\infty$ is sub-optimal, and we think that with some extra effort the results presented in this work can be adapted to a more classical (and more powerful) condition like $a_n/t_n\to\infty$, where $t_n$ holds for the number of final cells in the trees. In this case it is possible that not every tree would be consistent but the random forest would still converge to the regression function.

\section{Proof of \Cref{T1}}
\label{sec:Proofs}

First of all, we define the theoretical counterpart of the empirical CART criterion by
\begin{align*}
L^\star(A,d)=&\V[Y|\X\in A]-\V[Y|\X^{(h)}<z,\X\in A]\P[\X^{(h)}<z|\X\in A]\\
&\hphantom{{}\V[Y|\X\in A]}-\V[Y|\X^{(h)}\geq z,\X\in A]\P[\X^{(h)}\geq z|\X\in A].
\end{align*}
Analogously to the empirical case, we define the best theoretical cut $d^\star=(h^\star,z^\star)$ in $A$ as
\[d^\star=(h^\star,z^\star)\in\argmax_{\begin{tinymatrix}
d\in \mathcal{C}_A
\end{tinymatrix}} L^\star(A,d).\]
By the classical strong law of large numbers, $L_n(A,d)$ converges almost surely to \(L^\star(A,d)\) as $n$ tends to infinity, for all cuts $d\in\mathcal{C}_A$. This fact leads to the interpretation that the chosen cut at each step of the tree construction tends to decrease the variability of the sets of data points corresponding to the resulting nodes. This implies that the empirical mean in each cell tends to stabilize around a value that correspond to the conditional mean for the cell.
Applying basic algebra, it is easy to get an equivalent expression for the empirical CART criterion, given by
\begin{equation}
L_n(A,d)=\frac{N(A_L)N(A_R)}{N(A)N(A)}\parent{\bar{Y}_{A_L}-\bar{Y}_{A_R}}^2.\label{eq:CART}
\end{equation}
From \Cref{eq:CART}, we can get an alternative expression to the theoretical CART criterion, given by
\begin{align*}
L^\star(A,d)=&\P[\X^{(h)}<z|\X\in A]\P[\X^{(h)}\geq z|\X\in A]\\
&\times\Bigl(\E[Y|\X<z, \X\in A]-\E[Y|\X\geq z, \X\in A]\Bigr)^2.
\end{align*}
Remember that, in our procedure, the CART criterion takes ``imputed'' intervals which are updated after a cut is selected. From a theoretical point of view we do not longer have observations. Hence, we introduce the notions of the input assigned random variable $\widehat{\X}_{in}$ and output assigned random variable $\widehat{\X}_{out}$. The random variable $\widehat{\X}_{in}$ corresponds to a prior assignation whereas $\widehat{\X}_{out}$ corresponds to the assignations after the theoretical cut is performed. Furthermore, the binary assignations \(w\) are translated into probabilities.

Formally, let $\mathcal{W}$ be the collection of functions from $\R$ to $[0,1]$, and $\widehat{\X}_{in}=\parent{\widehat{\X}_{in}^{(1)},\ldots,\widehat{\X}_{in}^{(p)}}$ be the input distribution of the imputation, then $\widehat{\X}_{out}^{(h)}|\widehat{\X}_{in}$ is defined as
\[\widehat{\X}^{(h)}_{out}|\widehat{\X}_{in}\in A=\left\{\begin{array}{ll}
\X^{(h)}|\X\in A & \text{if } \M^{(h)}=0\\
\B^{(h)} & \text{if } \M^{(h)}=1\\
\end{array}, \right.\]
where
\begin{equation}
\B^{(h)}=\left\{\begin{array}{ll}
\parent{a^{(h)},z} & \text{if }\text{Ber}(w(Y))=1\\
\parent{z,b^{(h)}} & \text{if }\text{Ber}(w(Y))=0\\
\end{array} \right.,\quad w\in\mathcal{W}.\label{eq:B_hat_X}
\end{equation}
The imputation variable $w(Y)$ is the (random) probability that $\widehat{\X}_{out}^{(h)}<z$ conditionally to $\M^{(h)}=1$ and $Y$.
Note that $\widehat{\X}_{in}$ always belongs to a cell $A$, so the above definition of $\widehat{\X}^{(h)}_{out}|\widehat{\X}_{in}\in A$ is well defined.
We define the theoretical CART over a cut $d=(h,z)$ and a function $w\in\mathcal{W}$ as
\begin{align*}
L^\star(A,d,w)=&\V[Y|\widehat{\X}_{in}\in A]-\V[Y|\widehat{\X}_{out}^{(h)}<z,\widehat{\X}_{in}\in A]\P[\widehat{\X}_{out}^{(h)}<z|\widehat{\X}_{in}\in A]\\
&\hphantom{{}\V[Y|\widehat{\X}_{in}\in A]}-\V[Y|\widehat{\X}_{out}^{(h)}\geq z,\widehat{\X}_{in}\in A]\P[\widehat{\X}_{out}^{(h)}\geq z|\widehat{\X}_{in}\in A,]
\end{align*}
and the best cut and assignation $\parent{d^\star,w^\star}$ is selected by maximizing $L^\star\parent{A,d,w}$ over $\mathcal{C}_{A}\times\mathcal{W}$, that is
\[\parent{d^\star,w^\star}\in\argmax_{\begin{tinymatrix}
d\in \mathcal{C}_{A}\\
w\in\mathcal{W}
\end{tinymatrix}} L^\star\parent{A,d,w}.\]
We define, for any subset $A\subset \mathcal{X}$, the variation of $m$ within $A$ as \[\Delta(m,A)=\sup_{\x,\x'\in A}\abs{m(\x)-m(\x')}.\]
Furthermore, we denote by $A_{s(n)}(\X,\Theta)$ the final cell of the tree built with the random variable $\Theta$ that contains $\X$, where $s(n)$ is the number of cuts necessary to construct the cell in the tree.
After the last step of the tree construction, we end with a collection of imputed values that corresponds to the last imputed sample \(\parent{\widehat{\X}_{1,out},\dots,\widehat{\X}_{n,out}}\). Each of these vectors are non ambiguously assigned to a specific final cell of the tree.
In all that follows, when we write an imputation \(\widehat{\X}\) without any specification of \textit{in} or \textit{out}, we refer to the final imputation.
The proof of \Cref{T1} relies on \Cref{P1} below, which states that the variation of the regression function within a cell of a random forest is small for $n$ large enough and helps us to control the error of our predictor. In the sequel, we will abuse of the notation and will not write the dependence of the cells over $\X$ and $\Theta$.
\begin{proposition}\label{P1}
Assume that H\ref{H1} and H\ref{H2} hold. Then,
\[
\Delta(m, A_{s(n)})\to 0 \quad \text{almost surely.}
\]
\end{proposition}
For all $\x\in[0,1]^p$, we denote by  $L^\star\parent{A_{s(n)},d,w}$ the theoretical CART criterion over the cell $A_{s(n)}$ evaluated at a cut $d\in\mathcal{C}_{A_{s(n)}}$ and assignation $w\in\mathcal{W}$.
Let \((d^\star,w^\star)\) be the optimal couple (cut, assignation) of the cell \(A_{s(n)}\) for the theoretical criterion and let \((\widehat{d},\widehat{w})\) be the optimal couple (cut, assignation) of the cell \(A_{s(n)}\) for the empirical criterion.

\begin{lemma}\label{L1}
Assume that H\ref{H1} and H\ref{H2} are satisfied and fix $\x\in[0,1]^p$. Then for all $\rho,\xi>0$, there exists $N\in\mathbb{N}^\star$ such that, for all $n\geq N$ if $\P\corchet{L^\star\parent{A_{s(n)},d^\star,w^\star}\leq\xi}\geq 1-\rho$, then
\[
\Delta(m,A_{s(n)}(\x)) \to 0 \quad \text{almost surely.}
\]
\end{lemma}

\begin{lemma}\label{L2}
Assume that H\ref{H1} and H\ref{H2} are satisfied and fix $\x\in[0,1]^p$ . Then for all $\rho,\xi>0$, there exists $N\in\mathbb{N}^\star$ such that, for all $n\geq N$

\[\P\corchet{L_n\parent{A_{s(n)},\widehat{d},\widehat{w}}\leq \xi}\geq 1-\rho\]
\end{lemma}

\subsection{Proof of \Cref{L1}}

\begin{techlemma}\label{TL1}
Assume that H\ref{H1} and H\ref{H2} are satisfied and for a cell \(A\), $L^\star(A,d,\tilde{w}_{A,d})\equiv 0$ for all cuts $d=(h,z)\in \mathcal{C}_{A}$, where \[\tilde{w}_{A,d}=\P\corchet{a^{(h)}\leq \X^{(h)}<z|\X\in A,\M^{(h)}=1}.\]
Then, $m$ is constant on the cell $A$.
\end{techlemma}

\subsubsection*{Proof} Without loss of generality, we will assume that the cut $d=(1,z)$ is performed in the first direction so that $h=1$ and that the bounds of the cell $A$ are $a$ on the left and $b$ on the right. We omit the direction $h$ and simply note $X$ instead of $\X^{(1)}$. Note that if $\widehat{\X}_{out}$ is assigned to the left child node using $\tilde{w}_{A,d}$, which does not depend on $Y$, then $\widehat{\X}_{out}$ follows the same distribution than $\X$ and then our theoretical CART criterion is similar to the usual one for every cut $d=(h,z)\in\mathcal{C}_A$, that is,

\[L^\star\parent{A,d,\tilde{w}_{A,d}}=\P\corchet{a\leq X< z|\X\in A}\P\corchet{z\leq X\leq b|\X\in A}\parent{\mu_{A_{L}}-\mu_{A_{R}}}^2\]
where the notation $\mu_{A_{L}}$ (resp. $\mu_{A_{R}}$) holds for the conditional expected value of $Y$ given $\widehat{\X}_{in}\in A$ and $ a\leq X< z$ (resp. $ z\leq X\leq b$).
Because $\X$ is uniformly distributed over $[0,1]^p$,
\[\P\corchet{a\leq X< z|\X\in A}=\frac{z-a}{b-a}\]
and
\[\P\corchet{z\leq X\leq b|\X\in A}=\frac{b-z}{b-a}.\]
Next, we need to understand the distribution of $Y$ conditionally to $\widehat{\X}_{in}\in A$. This random variable is a mixture of the values of $Y$ such that $\X\in A$ and the ones that where assigned to the cell $A$ through the notion of the variable $\widehat{\X}_{in}$. Since the vector $(\widehat{\X}_{in},Y)$ has a density, we can give a precise meaning to the quantity
\begin{equation}
  \label{eq:def_m_tilde}
  \tilde{m}(\x)=\E\corchet{Y|\widehat{\X}_{in}=\x},
\end{equation}
for all $\x\in A$.
By assumption H\ref{H1}, we introduce the notation \(Y^{(j)}=m_j(X^{(j)})+\epsilon^{(j)}\), where \(\epsilon^{(j)}\sim \mathcal{N}(0,\sigma^2/p^2)\) so that we have \(Y\sim\sum_{j=1}^pY^{(j)}\).
Since under the condition $\X\in A$, the random variable $\X$ is also uniformly distributed on the cell $A$, we have that
\begin{align*}
\mu_{A_{L}}&=\E\corchet{\sum_{j=1}^p Y^{(j)}|a\leq X<z, \widehat{\X}_{in}\in A}\\
&=\sum_{j\ge 2}\E\corchet{Y^{(j)}|\widehat{\X}_{in}\in A}+\E\corchet{Y^{(1)}|a\leq X<z,\widehat{\X}_{in}\in A}\\
&=\tilde{K}+\frac{1}{z-a}C_a^z,
\end{align*}
where $\tilde{K}$ does not depend on $z$, $C_x^y=\int\limits_x^y\tilde{m}_1(t)dt$ and \\ \(\tilde{m}_1(t)=\E\corchet{Y^{(1)}|\widehat{\X}_{in}\in A,\ \X^{(1)}=t}\).
Analogously, we have that
\[\mu_{A_{R}}=\tilde{K}-\frac{1}{b-z}C_a^z+\frac{1}{b-z}C_a^b.\]
Therefore,
\begin{align*}
L^\star\parent{A,d,\tilde{w}_{A,d}}=&\parent{\frac{z-a}{b-a}}\parent{\frac{b-z}{b-a}}\parent{\frac{1}{z-a}C_a^z+\frac{1}{b-z}C_a^z-\frac{1}{b-z}C_a^b}^2\\\\
&=\frac{1}{(z-a)(b-z)}\parent{C_a^z-\frac{z-a}{b-a}C_a^b}^2.
\end{align*}
Since $L^\star\parent{A,d,\tilde{w}_{A,d}}=0$ by assumption, we obtain that for any $a\leq z\leq b$,
\[C_a^z = \frac{z-a}{b-a}C_a^b. \]
This proves that $z\mapsto C_a^z$ is linear in $z$ and thus, $\tilde{m}_1$ is constant on $\corchet{a,b}$. Since the law of \(Y_1\) is a mixture of distribution such that \(\M^{(1)}=0\) and \(\M^{(1)}=1\), we have that
\[
\tilde{m}_1(t)=(1-p^{(1)})m_1(t)+p^{(1)}\E\corchet{Y^{(1)}|\widehat{\X}_{in}\in A,\ \X^{(1)}=t,\ \M^{(1)}=1}
\]
where the second term does not depend on \(t\) (since the value is missing). This forces the function \(m_1\) to be constant on the interval \([a,b]\).
And, by additivity, the function $m$ is constant on $A$.

We need a second technical lemma that states that any cell within a cell for which each split gives a small value of $L^\star$ has also a uniformly small value of the CART criterion $L^\star$.

\begin{techlemma}
\label{TL_continuity_L}
Assume constructed the imputation \(\widehat{\X}_{in}\).
Let $\epsilon>0$, then there exists $\delta>0$ such that for all cell $A$ and cut $d\in\mathcal{C}_A$  satisfying $L^\star(A,d,\tilde{w}_{A,d})\le \delta$, we have
\[
L^\star(B,d,\tilde{w}_{A,d})\le \epsilon
\]
for all cell $B \subset A$, where both notions of CART criteria are constructed with the same vector of imputation \(\widehat{\X}_{in}\).
\end{techlemma}

\subsubsection*{Proof}
 We first prove the unidimensional case. Without loss of generality, we can assume that the cell $A$ is the interval $[0,1]$ and we first consider the cell $B$ of the form $[0,a]$ with $0\le a\le 1$. Following the lines of Technical Lemma \ref{TL1}, we have that
\[
L^\star(A,z,\tilde{w}_{A,d})=\frac{1}{z(1-z)}\left(C_0^z-C_0^1z\right)^{2},
\]
where for all $x$ and $y$, $C_x^y=\int_x^y \tilde{m}(t)dt$ and \(\tilde{m}\) is defined as in Equation \eqref{eq:def_m_tilde}. Evaluating the CART criterion at the value $z=a$, we have
\[
\frac{1}{a(1-a)}\left(C_0^a-C_0^1a\right)^{2}\le \delta.
\]
For the cell $B$, at a cut level $0<z<a$, we have that
\begin{align*}
    L^\star(B,z,\tilde{w}_{A,d})&=\frac{1}{z(a-z)}\left(C_0^z-C_0^a\frac{z}{a}\right)^{2}\\
    &=\frac{1}{z(a-z)}\left(C_0^z-C_0^1 z+C_0^1 z-C_0^a\frac{z}{a}\right)^{2}\\
    &\le \frac{2}{z(a-z)}\left(C_0^z-C_0^1 z\right)^2+\frac{2z}{a-z}\left(C_0^1 -\frac{C_0^a}{a}\right)^{2}\\
    &\le 2\frac{1-z}{a-z}L^\star(A,z,\tilde{w}) + \frac{2z(1-a)}{a(a-z)}\delta\\
    &\le 2\frac{a(1-z)+z(1-a)}{a(a-z)}\delta \le  \frac{2}{a(a-z)}\delta.
\end{align*}
But since the function $z\mapsto C_0^z$ is differentiable, we have that $C_0^z=C_0^a-(a-z)\tilde{m}(a)+o(a-z)$ when $z\to a$ by a Taylor expansion. This shows that $L^\star(B,z,\tilde{w}_{A,d}) \to 0$ when $z\to a$ and then there exists a $\delta_0>0$ such that $L^\star(B,z,\tilde{w}_{A,d})\le \epsilon$ when $(a-z)\le \delta_0$. Then, using the previous inequalities for $0<z<a-\delta_0$, we have that for $\delta=\epsilon\delta_0 a/2$, $L^\star(B,z,\tilde{w}_{A,d})\le \epsilon$. In the same way, we generalize the previous arguments for $A=[0,1]^p$ and for a specific type of cell
\[
B=\{\x\in A: \x^{(1)}\le a\}
\]
with $0<a<1$. The values $C_y^z$ are then replaced by $\int_y^z\int_0^1\dots\int_0^1\tilde{m}(t)dt_1dt_2\dots dt_p$. The result can be repeated for the case $B=\{\x\in A: \x^{(1)}\ge a\}$.

For the general case of $B\subset A$ we see that any $B$ can be obtained by a finite sequence of $B=B_k\subset B_{k-1}\subset \dots B_1\subset A$ of subset constructed by the scheme described above. This finishes the proof.

\subsubsection*{Proof of \Cref{L1}}

We will show that $\Delta(m,A_{s(n)}(\x)) \to 0$ a.s. by contradiction.
We assume that with positive probability, there exists a positive constant $c>0$ and a sub-sequence $\phi(n)$ of cells  $A_{s(\phi(n))}$ such that $\Delta(m,A_{s(\phi(n))}(\x))> c$. This means that in each set $A_{s(\phi(n))}(\x)$ one can find a pair of elements $(x_n,y_n)$ such that
\[
|m(x_n)-m(y_n)|>c.
\]
The sequences $(x_n)_n$ and $(y_n)_n$ belong to the compact set $[0,1]^p$ so one can extract a sub-sequence $\psi(n)$ such that $(x_n)_n$ and $(y_n)_n$ converge respectively to two points $x$ and $y$. By continuity of $m$, we have that
\[
|m(x)-m(y)|>c.
\]
At the cost of taking an $x'$ and $y'$ close to $x$ and $y$, satisfying
\[
|m(x')-m(y')|>\frac{c}{2},
\]
and such that, for $n$ large enough, all the cells $A_{s(\phi(n))}$ contain the pair of points $x'$ and $y'$.
By hypothesis, we know that
\[
\sup_{d\in\mathcal{C}_{A_{s(\phi\circ\psi(n))}}}L^\star(A_{s(\phi\circ\psi(n))},d,\tilde{w})\to 0 \quad \text{in probability},
\]
where we just wrote $\tilde{w}$ for the choice of the assignation given in Technical Lemma \ref{TL1}.
So one can extract a subsequence \(\chi(n)\) such that the \(\sup_{d}L^\star(A_{s(\phi\circ \psi\circ \chi(n))},d,\tilde{w})\) converges to 0 almost surely. For simplicity, we keep denoting \(n\) for the sub-sequence \(\phi\circ \psi\circ \chi(n)\) in the following part of the proof.
 Lastly, we define the sequence of cells $(C_i)_{i\ge 1}$ such that
\[
C_i=\bigcap_{k=1}^{i}A_{s(k)}.
\]
These cells form a non-increasing sequence for the inclusion order and for each \(i\), \(C_i\subset A_{s(i)}\). The cells $C_i$ inherit the same vector of imputation as for $A_{s(i)}$. We can use Technical Lemma \ref{TL_continuity_L} to get that \(\sup_{d\in\mathcal{C}_{A_{s(n)}}}L^\star(C_n,d,\tilde{w}) \to 0\) a.s. Since the sequence of cells \((C_i)_i\) is a non-increasing sequence, there exists a cell, denoted \(C_{\infty}\) that is the limit of the cells \(C_i\) when \(i\to \infty\) in the sense
\begin{equation*}
  C_{\infty}=\bigcap_{i\ge 1}C_i.
\end{equation*}
To see that, one can write \(C_i=\prod_h [a_i^{(h)},b_i^{(h)}]\) and take the limits of the sequences \((a_i^{(h)})_i\) and \((b_i^{(h)})_i\).
The objective is to show that there is no other possibility than having \(L^\star(C_{\infty},d,\tilde{w})=0\) for every cut \(d\).
% This is done by repeating the ideas behind the proof of Lemma 1 in \cite{scornet2015consistency}.
For a cell \(A\), the CART criterion \(L^\star(A,d,\tilde{w}_{A,d})\) has a continuous behavior with respect to \(A\). Indeed, we see that
\begin{align*}
L^\star\parent{A,d,\tilde{w}_{A,d}}=&\P\corchet{a^{(h)}\leq \X^{(h)}< z|\X\in A}\P\corchet{z\leq \X^{(h)}\leq b^{(h)}|\X\in A}\parent{\widehat{\mu}_{A_L}-\widehat{\mu}_{A_R}}^2,
\end{align*}
which is a product of three terms that are uniformly continuous in \(A\) (since \(m\) is a continuous function) with respect to the distance between cells given by \(\Delta(A,B)=\max_{h}|x^{(h)}_A-x_B^{(h)}|+\max_{h}|y^{(h)}_A-y_B^{(h)}|\) where \(x_A^{(h)}\) and \(y_A^{(h)}\) are defined as \(A=\bigcap_{h}[x_A^{(h)},y_A^{(h)}]\).
Since \(C_i\to C_{\infty}\) for the distance \(\Delta\), we have that for all \(\epsilon>0\), for all \(i\) large enough, for all cut \(d\) of the final cell \(C_{\infty}\),
\[
|L^\star\parent{C_i,d,\tilde{w}}-L^\star\parent{C_{\infty},d,\tilde{w}}| \le \epsilon.
\]
But since the \(L^\star\parent{C_i,d,\tilde{w}}\) converges almost surely uniformly in \(d\) to 0, and \(\epsilon\) is arbitrary,  for every cut \(d\) inside the valid cuts of \(C_{\infty}\), we have that \(L^\star(C_{\infty},d,\tilde{w})=0\). But then by Technical Lemma \ref{TL1}, the function \(m\) has to be constant in the cell \(C_{\infty}\). But the points \(x'\) and \(y'\) do belong to the cell \(C_{\infty}\) since they belong to each of the cells in the intersection. So \(m(x')=m(y')\) which contradicts the fact that \(|m(x')-m(y')|>c/2\). This proves that
\[
\Delta(m,A_{s(n)}(\x)) \to 0 \quad \text{almost surely.}
\]

\subsection{Proof of \Cref{L2}}

Remember that $A_{s(n)}$ denotes the final cell of the tree built with the random variable $\Theta$ that contains $\X$, where $s(n)$ is the number of cuts necessary to construct the cell. Similarly, $A_k$ is the same cell but where only the first $k$ cuts have been performed.

\begin{techlemma}
Assume that H\ref{H1} and H\ref{H2} hold and fix $\x\in[0,1]^p$.  For all $\rho,\xi>0$ there exists $N\in\mathbb{N}^\star$ such that for all $n\geq N$ there exists $k_0(n)\in\mathbb{N}^\star$ such that
\[\P\corchet{\bigl|L_n\parent{A_{s(n)},\widehat{d},\widehat{w}}-L_n\parent{A_k,\widehat{d},\widehat{w}}\bigr|\leq \xi}\geq 1-\rho,\]
for all  $k\geq k_0(n)$.
\label{TL2}
\end{techlemma}

\subsubsection*{Proof}
Fix $\alpha,\rho>0$ and consider the following standard inequality on a Gaussian tail
\[\P[\varepsilon_1\ge \alpha]\le \frac{\sigma^2}{t\sqrt{2\pi}}\exp\parent{-\frac{t^2}{2\sigma^2}}\]
then, simple calculations show that, for all $n\in\N^\star$

\begin{equation}
\P\llavs{\left|\sum_{i=1}^n\varepsilon_i\right|\geq n\alpha}\le \frac{\sigma}{\alpha\sqrt{n}}\exp\llavs{-\frac{\alpha^2n}{2\sigma^2}}.\label{eq:random_mean_epsilon}
\end{equation}

Note that there are at most $n(n+1)/2$ sets of the form $\{i:\X_i^{(h)}\in[a_n,b_n],\M_i^{(h)}=0\}$ for $0\le a_n<b_n\le 1$. On the other hand, let $(Y_{(1)},\ldots,Y_{(n)})$ be the ordered vector of $Y$, since missing observations are assigned to the cell using $Y$ and maximizing the CART criterion, this implies that close values of $Y$ must be assigned to the same cell, thus once again note that there are at most $n(n+1)/2$ sets of the form $\{i:Y_{(i)}\in[a_n,b_n],\M_i^{(h)}=1\}$ for $0\le a_n<b_n\le 1$.

We deduce from \Cref{eq:random_mean_epsilon} and the union bound, that there exists $N_1\in\N^\star$ such that, with probability at least $1-\rho$, for all $n\ge N_1$ and all $0\le a_n<b_n\le 1$ satisfying $\widehat{N}\parent{\prod_{h=1}^p [a_n^{(h)},b_n^{(h)}]}\ge q_n$,

\begin{equation}
\Biggl|\frac{1}{\widehat{N}\parent{\prod_{h=1}^p [a_n^{(h)},b_n^{(h)}]}}\sum_{i=1}^n\varepsilon_i\indicadora{\widehat{\X}_i\in A}{}\Biggr|\le \frac{\sigma n^{4p}}{\alpha\sqrt{n}}\exp\llavs{-\frac{\alpha^2n}{2\sigma^2}} \leq\alpha.\label{eq:lemma2_3}
\end{equation}
Furthermore, making use of the same ideas and applying the inequality $\P[\chi^2(n)\ge 5n]\le e^{-n}$ (for interested readers, see \citet{laurent2000adaptive}), there exists $N_2\in\mathbb{N}^\star$ such that, with probability at least $1-\rho$ for all $n\ge N_2$ and all $0\le a_n<b_n\le 1$ satisfying $\widehat{N}\parent{\prod_{h=1}^p [a_n^{(h)},b_n^{(h)}]}\ge q_n$
\begin{equation}
\frac{1}{\widehat{N}\parent{\prod_{h=1}^p [a_n^{(h)},b_n^{(h)}]}}\sum_{i=1}^n\varepsilon_i^2\indicadora{\widehat{\X}_i\in A}{}\leq\tilde{\sigma}^2,\label{eq:lemma2_4}
\end{equation}
where $\tilde{\sigma}^2$ is a positive constant, depending only on $\rho$.
Since $(A_{k})_k$ is a decreasing sequence of compact sets, for every $\xi>0$ there exists $k_0$ such that, for all $k\geq k_0$
\begin{equation}
\max\parent{\norm{\mathbf{a}_{k}-\mathbf{a}_{s(n)}}_\infty, \norm{\mathbf{b}_{k}-\mathbf{b}_{s(n)}}_\infty}\leq \xi,\label{eq:lemma1_maxdelta}
\end{equation}
where $\mathbf{a}_{k}=\parent{a_{k}^{(1)},\ldots,a_{k}^{(p)}}\in[0,1]^p$ and $\mathbf{b}_{k}=\parent{b_{k}^{(1)},\ldots,b_{k}^{(p)}}\in[0,1]^p$ such that
\[A_{k}=\prod_{j=1}^p\corchet{a_{k}^{(j)},b_{k}^{(j)}}.\]
We define $\mathbf{a}_{s(n)}$ and $\mathbf{b}_{s(n)}$ analogously for the cell $A_{s(n)}$.
Since the variable $\X$ is uniformly distributed in the hypercube $[0,1]^p$ and the missing entries follow the MCAR mechanism we have that, for any cell $A$, $\widehat{N}_{obs}(A)/\widehat{N}_{obs}([0,1]^p)\to \text{Vol}(A)$ almost surely as $n\to \infty$. Furthermore,
for any $k\geq k_0(n)$, Equation \eqref{eq:lemma1_maxdelta} implies that $\text{Vol}(A_{k})\leq (1+2\xi)^p\text{Vol}(A_{s(n)})$. Then there exists $N_3\in \N^\star$ such that for all $n\geq N_3$, $\widehat{N}_{obs}(A_{k})\leq (1+3\xi)^p \widehat{N}_{obs}(A_{s(n)})$, and
\begin{equation}
\widehat{N}_{obs}(A_{k}\setminus A_{s(n)})\leq\xi' \widehat{N}_{obs}(A_{s(n)})\leq\xi' \widehat{N}_{obs}(A_{k}),\label{eq:NN}
\end{equation}
where $\xi'=(1+3\xi)^p-1$.
On the other hand, for any fixed $n$, the quantity  $u_k=\widehat{N}_{miss}(A_{k})-\widehat{N}_{miss}(A_{s(n)})$ converges to 0 almost surely as $k\to \infty$.
Now taking $n\ge N_3$ fixed, we have that there exists a $k_1(n)$ such that for all $k\ge k_1(n)$, with probability at least $1-\rho$, we have
\begin{align}
\widehat{N}(A_{k}\backslash A_{s(n)})&=(\widehat{N}_{obs}(A_{k})+\widehat{N}_{miss}(A_{k}))-(\widehat{N}_{obs}(A_{s(n)})+\widehat{N}_{miss}(A_{s(n)}))\nonumber\\
&=\widehat{N}_{obs}(A_{k}\setminus A_{s(n)})+u_k\nonumber\\
&\leq 2\xi' \widehat{N}_{obs}(A_{s(n)}),\label{eq:N_Vol}
\end{align}
where we used Equation \eqref{eq:NN} and the fact that convergence almost sure implies convergence in probability.
For the rest of the proof, we take $k\geq \max\{k_0(n),k_1(n)\}$ and assume that Equations \eqref{eq:lemma2_3},\eqref{eq:lemma2_4} and \eqref{eq:N_Vol} are satisfied, which occurs with probability at least $1-3\rho$ for every $n>N$ with $N=\max\{N_1,N_2,N_3\}$.
Note that $A_{k}\setminus A_{s(n)}$ is either a final node, contains at least one final node or is empty (in which case the result holds trivially). Since each final node contains at least $q_n$ points then for $q_n$ sufficiently large, using \cref{eq:lemma2_3} we conclude that $|\widehat{Y}_{A_{k}}|\leq \norm{m}_\infty+\alpha$, $|\widehat{Y}_{A_{s(n)}}|\leq \norm{m}_\infty+\alpha$ and $|\widehat{Y}_{A_{k}\setminus A_{s(n)}}|\leq \norm{m}_\infty+\alpha$.
We use the following decomposition
\[\bigl|L_n\parent{A_{k},\widehat{d},\widehat{w}}-L_n\parent{A_{s(n)},\widehat{d},\widehat{w}}\bigr|\leq K_0+K_L+K_R,\]
where the three terms $K_0$, $K_L$ and $K_R$ are given by
\[K_0 =\Big|\frac{1}{\widehat{N}(A_{k})}\sum_{i=1}^n (Y_i-\widehat{Y}_{A_{k}})^2\indicadora{\widehat{\X}_i\in A_{k}}{}-\frac{1}{\widehat{N}(A_{s(n)})}\sum_{i=1}^n (Y_i-\widehat{Y}_{A_{s(n)}})^2\indicadora{\widehat{\X}_i\in A_{s(n)}}{}\Big|,\]
\begin{align*}
K_L =&\Big|\frac{1}{\widehat{N}(A_{k})}\sum_{i=1}^n (Y_i-\widehat{Y}_{A_{L,k}})^2\indicadora{\widehat{\X}_i\in A_{k},\widehat{\X}^{(\widehat{h})}_i<\widehat{z}}{}\\
&\hphantom{\Big|\frac{1}{\widehat{N}(A_{k})}\sum_{i=1}^n (Y_i-\widehat{Y}_{A_{L,k}})^2}-\frac{1}{\widehat{N}(A_{s(n)})}\sum_{i=1}^n (Y_i-\widehat{Y}_{A_{L,s(n)}})^2\indicadora{\widehat{\X}_i\in A_{s(n)},\widehat{\X}^{(\widehat{h})}_i<\widehat{z}}{}\Big|
\end{align*}
and
\begin{align*}
K_R =&\Big|\frac{1}{\widehat{N}(A_{k})}\sum_{i=1}^n (Y_i-\widehat{Y}_{A_{R,k}})^2\indicadora{\widehat{\X}_i\in A_{k},\widehat{\X}^{(\widehat{h})}_i\geq\widehat{z}}{}\\
&\hphantom{\Big|\frac{1}{\widehat{N}(A_{k})}\sum_{i=1}^n (Y_i-\widehat{Y}_{A_{R,k}})^2}-\frac{1}{\widehat{N}(A_{s(n)})}\sum_{i=1}^n (Y_i-\widehat{Y}_{A_{R,s(n)}})^2\indicadora{\widehat{\X}_i\in A_{s(n)},\widehat{\X}^{(\widehat{h})}_i\geq\widehat{z}}{}\Big|.
\end{align*}
We first consider the term $K_0$ that can be upper bounded once again by using a similar split in $K_0\leq K_{0,1}+K_{0,2}+K_{0,3}$ where
\[K_{0,1}=\Big|\frac{1}{\widehat{N}(A_{k})}\sum_{i=1}^n (Y_i-\widehat{Y}_{A_{k}})^2\indicadora{\widehat{\X}_i\in A_{s(n)}}{}-\frac{1}{\widehat{N}(A_{k})}\sum_{i=1}^n (Y_i-\widehat{Y}_{A_{s(n)}})^2\indicadora{\widehat{\X}_i\in A_{s(n)}}{}\Big|,\]
\[K_{0,2}=\Big|\frac{1}{\widehat{N}(A_{k})}\sum_{i=1}^n (Y_i-\widehat{Y}_{A_{s(n)}})^2\indicadora{\widehat{\X}_i\in A_{s(n)}}{}-\frac{1}{\widehat{N}(A_{s(n)})}\sum_{i=1}^n (Y_i-\widehat{Y}_{A_{s(n)}})^2\indicadora{\widehat{\X}_i\in A_{s(n)}}{}\Big|\]
and
\[K_{0,3}=\Big|\frac{1}{\widehat{N}(A_{k})}\sum_{i=1}^n (Y_i-\widehat{Y}_{A_{k}})^2\indicadora{\widehat{\X}_i\in A_{k}\setminus A_{s(n)}}{}\Big|.\]
For $K_{0,1}$, observe that
\begin{align*}
|\widehat{Y}_{A_{k}}-\widehat{Y}_{A_{s(n)}}|&
=\Big|\frac{1}{\widehat{N}(A_{k})}\sum_{i=1}^n Y_i\indicadora{\widehat{\X}_i\in A_{k}\setminus A_{s(n)}}{}+\frac{1}{\widehat{N}(A_{k})}\sum_{i=1}^n Y_i\indicadora{\widehat{\X}_i\in A_{s(n)}}{}-\widehat{Y}_{A_{s(n)}}\Big|\\
&\leq \frac{\widehat{N}(A_{k}\setminus A_{s(n)})}{\widehat{N}(A_{k})}|\widehat{Y}_{A_{k}\setminus A_{s(n)}}-\widehat{Y}_{A_{s(n)}}|\\
&\leq 4\xi'(\norm{m}_\infty+\alpha).
\end{align*}
Hence,
\begin{align*}
K_{0,1}&
\leq \frac{2}{\widehat{N}(A_{k})}|\widehat{Y}_{A_{s(n)}}-\widehat{Y}_{A_{k}}|\left|\sum_{i=1}^n \parent{Y_i+\frac{\widehat{Y}_{A_{s(n)}}+\widehat{Y}_{A_{k}}}{2}}\indicadora{\widehat{\X}_i\in A_{s(n)}}{}\right|\\
&\leq\frac{8\xi'(\norm{m}_\infty+\alpha)}{\widehat{N}(A_{k})}\corchet{\left|\widehat{N}(A_{s(n)})\widehat{Y}_{A_{s(n)}}\right|+\left|\frac{\widehat{Y}_{A_{s(n)}}+\widehat{Y}_{A_{k}}}{2} \widehat{N}(A_{s(n)})\right|}\\
&\leq 16\xi'(\norm{m}_\infty+\alpha)^2.
\end{align*}
For the term $K_{0,2}$, with the help of Equation \eqref{eq:lemma2_4} observe that
\begin{align*}
K_{0,2}&
\leq \frac{\widehat{N}(A_k\setminus A_{s(n)})}{\widehat{N}(A_k)}\Big|\frac{1}{\widehat{N}(A_{s(n)})}\sum_{i=1}^n (Y_i-\widehat{Y}_{A_{s(n)}})^2\indicadora{\widehat{\X}_i\in A_{s(n)}}{}\Big|\\
&\leq 2\xi'\Big|\frac{1}{\widehat{N}(A_{s(n)})}\sum_{i=1}^n Y_i^2\indicadora{\widehat{\X}_i\in A_{s(n)}}{}+\widehat{Y}^2_{A_{s(n)}}\Big|\\
&\leq 2\xi'\Bigg[(\norm{m}_\infty+\alpha)^2+\frac{1}{\widehat{N}(A_{s(n)})}\sum_{i=1}^n m^2(\X_i)\indicadora{\widehat{\X}_i\in A_{s(n)}}{}\\
&\hphantom{\leq \xi\Bigg[}+\Big|\frac{2}{\widehat{N}(A_{s(n)})}\sum_{i=1}^n m(\X_i)\varepsilon_i\indicadora{\widehat{\X}_i\in A_{s(n)}}{}\Big|+\frac{1}{\widehat{N}(A_{s(n)})}\sum_{i=1}^n \varepsilon_i^2\indicadora{\widehat{\X}_i\in A_{s(n)}}{}\Bigg]\\
&\leq 2\xi'\corchet{(\norm{m}_\infty+\alpha)^2+\norm{m}_\infty^2+2\norm{m}_\infty\alpha+\tilde{\sigma}^2}.
\end{align*}
Regarding $K_{0,3}$, observe that
\begin{align*}
K_{0,3}&
\leq 2\xi'\Big|\frac{1}{\widehat{N}(A_k\setminus A_{s(n)})}\sum_{i=1}^n(Y_i^2+2Y_i\widehat{Y}_{A_k}+\widehat{Y}_{A_k}^2)\indicadora{\widehat{\X}_i\in A_k\setminus A_{s(n)}}{}\Big|\\
&\leq 2\xi'\Bigg[\Big|\frac{1}{\widehat{N}(A_k\setminus A_{s(n)})}\sum_{i=1}^n (m(\X_i)+\varepsilon_i)^2\indicadora{\widehat{\X}_i\in A_k\setminus A_{s(n)}}{}\Big|\\
&\hphantom{\leq \xi\Bigg[}+\frac{2}{\widehat{N}(A_k\setminus A_{s(n)})}(\norm{m}_\infty+\alpha)\Big|\sum_{i=1}^n (m(\X_i)+\varepsilon_i)\indicadora{\widehat{\X}_i\in A_k\setminus A_{s(n)}}{}\Big|\\
&\hphantom{\leq \xi\Bigg[}+(\norm{m}_\infty+\alpha)^2\Bigg]\\
&\leq 2\xi'\Bigg[\norm{m}_\infty^2+2\norm{m}_\infty\Big|\frac{1}{\widehat{N}(A_k\setminus A_{s(n)})}\sum_{i=1}^n\varepsilon_i\indicadora{\widehat{\X}_i\in A_k\setminus A_{s(n)}}{}\Big|\\
&\hphantom{\leq \xi\Bigg[}+\frac{1}{\widehat{N}(A_k\setminus A_{s(n)})}\sum_{i=1}^n\varepsilon_i^2\indicadora{\widehat{\X}_i\in A_k\setminus A_{s(n)}}{}+2(\norm{m}_\infty+\alpha)^2+(\norm{m}_\infty+\alpha)^2\Bigg]\\
&\leq 2\xi'[3(\norm{m}_\infty+\alpha)^2+\norm{m}_\infty^2+2\norm{m}_\infty\alpha+\tilde{\sigma}^2].
\end{align*}
Therefore, there exists a universal constant $C>0$ such that $K_0\leq C\xi$ and with similar arguments we can show that $K_L\leq C\xi$ and $K_R\leq C\xi$. Which concludes that with probability at least $1-3\rho$,
\[\bigl|L_n\parent{A_k,\widehat{d},\widehat{w}}-L_n\parent{A_{s(n)},\widehat{d},\widehat{w}}\bigr|\leq 3C\xi.\]

\subsubsection*{Proof of \Cref{L2}}

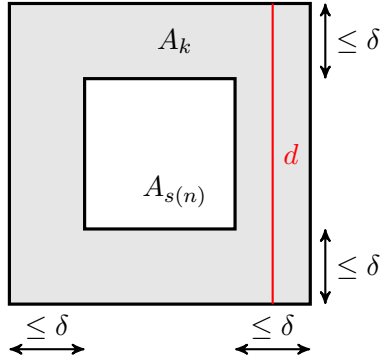
\begin{figure}[ht]
\begin{center}
\begin{minipage}{0.45\textwidth}
\begin{center}
	\begin{tikzpicture}[scale=4, >=stealth']
	% Aqui esta la celda A_k
	\draw[Import_line, fill = black!10!white] (0,0) rectangle (1,1);
	\node at (0.55,0.875) {$A_k$};

	% Aqui esta la celda A_infty
	\draw[Import_line, fill = white] (0.25,0.25) rectangle (0.75,0.75);
	\node at (0.55,0.375) {$A_{s(n)}$};

	\draw[thick, <->] (0,-0.15) -- node[above] {$\leq \delta$} (0.25,-0.15);
	\draw[thick, <->] (0.75,-0.15) -- node[above] {$\leq \delta$} (1,-0.15);

	\draw[thick, <->] (1.05,0) -- node[right] {$\leq \delta$} (1.05,0.25);
	\draw[thick, <->] (1.05,0.75) -- node[right] {$\leq \delta$} (1.05,1);

	% Aqui esta posible y no posible corte
	\draw[thick, color = red] (0.875,0) -- node[right] {$d$} (0.875,1);
	%\draw[thick, color = red] (0.3,0) -- node[right] {$d'$} (0.3,1);
	\end{tikzpicture}
\caption{All feasible cuts in $A_k$ must be performed in the $A_k\setminus A_{s(n)}$, like the cut $d$ in the figure, otherwise the cut would split $A_{s(n)}$.}
\label{fig:lemma1_2}
\end{center}
\end{minipage}
\end{center}
\end{figure}

Assume that H\ref{H1} and H\ref{H2} are satisfied, fix $\x\in[0,1]^p$ and $\Theta$. Then, let us show by  contradiction that for all $\xi>0$, there exists $N\in\mathbb{N}^\star$ such that, with probability at least $1-\rho$ for all $n\geq N$
\[L_n\parent{A_{s(n)},\widehat{d},\widehat{w}}\leq \xi.\]
So, assume that there exists $c>0$, $0<p_0<1$ and a sub-sequence $\phi(n)$ such that
\[L_{\phi(n)}\parent{A_{s(\phi(n))},\widehat{d}_{s(\phi(n))},\widehat{w}_{s(\phi(n))}}>c,\]
with probability at least $p_0$. To keep the notation simple, we omit to write $\phi(n)$ and still write $n$ for the indexes of the sub-sequence.
Additionally, assume that $k$ is sufficiently large so the conclusion of Technical Lemma \ref{TL2} is satisfied, hence
\begin{equation}
\bigl|L_n\parent{A_k,\widehat{d},\widehat{w}}-L_n\parent{A_{s(n)},\widehat{d},\widehat{w}}\bigr|\leq \xi\label{eq:lemma1_supL_0}
\end{equation}
and \Cref{eq:lemma1_maxdelta} is satisfied.
Note that all the feasible cuts $d$ in $A_k$ must be performed in $A_k\setminus A_{s(n)}$, otherwise $d$ would split $A_{s(n)}$ and $(A_k)_k$ would not converge to $A_{s(n)}$ (see \Cref{fig:lemma1_2} for an illustration in $p=2$). From here we conclude that
\[L_n\parent{A_k,\widehat{d},\widehat{w}} \leq \sup_{\begin{tinymatrix}
d\in \mathcal{C}_{A_k}\cap\mathcal{C}_{A_{s(n)}}\\
w\in\mathcal{W}_{A_k}
\end{tinymatrix}} L_n\parent{A_k,d,w} \leq \sup_{\begin{tinymatrix}
d\in \mathcal{C}_{A_k}\\
w\in\mathcal{W}_{A_k}
\end{tinymatrix}} L_n\parent{A_k,d,w}\leq \xi.\]

On the other hand, from \Cref{eq:lemma1_supL_0}, with probability at least $p_0$, we have
\[c-\xi\leq L_n\parent{A_{s(n)},\widehat{d},\widehat{w}}-\xi^2\leq L_n\parent{A_k,\widehat{d},\widehat{w}}.\]
Hence, we have, with probability at least $p_0$,
\[c-C\xi\leq L_n\parent{A_k,\widehat{d},\widehat{w}}\leq \xi\] which is absurd, since $c>0$ is fixed and $\xi$ is arbitrarily small. Thus the result follows.

\subsection{Proof of \Cref{P1}}
For a cell $A$, fix a cut $d\in\mathcal{C}_A$ and consider a function $w\in\mathcal{W}$, we need to define $L_n(A,d,w)$. This is done according to the following procedure, first create a random vector $W$ of dimension $\widehat{N}_{miss}^{(h)}(A)=\text{Card}(\imiss_{A,miss}^{(h)})$, where $W_k=\text{Ber}(w(Y_{j_k}))$ for $j_k\in\imiss_{A,miss}^{(h)}$, then assign the observations $\widehat{\X}_{j_k}$ to the child nodes according to the random vector $W$. Once we have assigned the observations to the child nodes, we evaluate the empirical CART criterion $L_n$ considering these assignations. Note that in this case $L_n$ is a random variable and the assignations are independent to each other so $L_n(A,d,w)$ is a sum of independent random variables with the same distribution. Hence, by the strong law of large numbers $L_n(A,d,w)\to L^\star(A,d,w)$ almost surely as $n\to\infty$ for all cuts $d\in\mathcal{C}_A$ and all functions $w\in\mathcal{W}$.

We prove the almost sure convergence of $\Delta(m,A_{s(n)})$ towards 0 by showing that the theoretical CART criterion of the sequence $(A_{s(n)})_n$ tends to 0 and use of \Cref{L1,L2}.
Note that
\begin{align*}
&L^\star\parent{A_{s(n)},d^\star,w^\star}-L_n\parent{A_{s(n)},\widehat{d},\widehat{w}}\\&=L^\star\parent{A_{s(n)},d^\star,w^\star}-L_n\parent{A_{s(n)},d^\star,w^\star}\\&\hphantom{L^\star\parent{A_{s(n)},d^\star,w^\star}}+L_n\parent{A_{s(n)},d^\star,w^\star}-L_n\parent{A_{s(n)},\widehat{d},\widehat{w}}\\
&\leq L^\star\parent{A_{s(n)},d^\star,w^\star}-L_n\parent{A_{s(n)},d^\star,w^\star}.
\end{align*}
Where the last inequality comes from noting that $L_n\parent{A_{s(n)},\widehat{d},\widehat{w}}\geq L_n\parent{A_{s(n)},d,w}$ for all cut $d\in\mathcal{C}_{A_{s(n)}}$ and assignation $w\in\mathcal{W}^{(\widehat{h})}_{A_{s(n)}}$, where $\widehat{d}=(\widehat{h},\widehat{z})$.
As discussed above, by strong law of large numbers \[L^\star\parent{A_{s(n)},d^\star,w^\star}-L_n\parent{A_{s(n)},d^\star,w^\star}\to 0,\] almost surely.
Fix $\xi,\rho>0$, for $n$ sufficiently large, we have
\[L^\star\parent{A_{s(n)},d^\star,w^\star}-L_n\parent{A_{s(n)},\widehat{d},\widehat{w}}\leq \xi\quad\text{almost surely.}\]
On the other hand, by \Cref{L2}, there exists $N_1$ such that for all $n\geq N_1$, with probability at least $1-\rho$
\[L_n\parent{A_{s(n)},\widehat{d},\widehat{w}}\leq C\xi.\]
Hence, with the same probability,
\[L^\star\parent{A_{s(n)},d^\star,w^\star}\leq\xi.\]
And by \Cref{L1}, we conclude that
\[\Delta(m,A_{s(n)})\xrightarrow{a.s.}0.\]

\subsection{Proof of \Cref{T1}}
Let \[m_n(\X)=\frac{1}{\widehat{N}(A_{s(n)}(\X))}\sum_{i=1}^n Y_i\indicadora{\widehat{\X}_i\in A_{s(n)}(\X)}{}\]
be our tree estimator. Moreover, let us define two other estimators, the first one takes our partition of the cells $A_{s(n)}$ built up using the imputed variables but considers the local means using the complete (unseen) observations $\X_i$, denoted as
\[m_n'(\X)=\frac{1}{N(A_{s(n)}(\X))}\sum_{i=1}^n Y_i\indicadora{\X_i\in A_{s(n)}(\X)}{},\]
while the second one takes $m(\X_i)$ for the prediction and the complete observations $\X_i$, denoted as
\[m_n''(\X)=\frac{1}{N(A_{s(n)}(\X))}\sum_{i=1}^n m(\X_i)\indicadora{\X_i\in A_{s(n)}(\X)}{}.\]
By \Cref{eq:lemma2_3}, we know that for all $\alpha,\xi>0$ there exists $N\in\N^\star$, such that for all $n\ge N$,
\begin{align*}
&\P\corchet{|m_n'(\X)-m_n''(\X)|\ge \alpha}\\
&=\P\Biggl[\Biggl|\frac{1}{N(A_{s(n)}(\X))}\sum_{i=1}^n \Bigl(Y_i-m(\X_i)\Bigr)\indicadora{\X_i\in A_{s(n)}(\X)}{}\Biggr|\ge \alpha\Biggr]\\
&=\P\Biggl[\Biggl|\frac{1}{N(A_{s(n)}(\X))}\sum_{i=1}^n \varepsilon_i\indicadora{\X_i\in A_{s(n)}(\X)}{}\Biggr|\Bigr]\\
&\le \xi.
\end{align*}
On the other hand, note that $m_n''(\X)=\sum_{i=1}^n W_{n,i}(\X)m(\X_i)$, where
\[W_{n,i}(\X)=\frac{1}{N(A_{s(n)}(\X))}\indicadora{\X_i\in A_{s(n)}(\X)}{}.\]
Then,
\begin{align*}
\E[m_n''(\X)-m(\X)]^2& = \E\corchet{\sum_{i=1}^n W_{n,i}(\X)m(\X_i)-m(\X)}^2\\
&=\E\Biggl[\sum_{i=1}^n\sqrt{W_{n,i}(\X)}\sqrt{W_{n,i}(\X)}(m(\X_i)-m(\X))\Biggr]^2\\
&\text{(Applying Cauchy-Schwartz's inequality)}\\
&\le\E\corchet{\sum_{i=1}^n W_{n,i}(\X)\sum_{i=1}^n W_{n,i}(\X)(m(\X_i)-m(\X))^2}\\
&=\E\corchet{\sum_{i=1}^n W_{n,i}(\X)(m(\X_i)-m(\X))^2\indicadora{\X_i\in A_{s(n)}(\X)}{}}.
\end{align*}
Note that $(m(\X_i)-m(\X))^2\indicadora{\X_i\in A_{s(n)}(\X)}{}\le \Delta(m,A_{s(n)}(\X))^2$, hence
\[\E[m_n''(\X)-m(\X)]^2\leq \E\corchet{\Delta(m,A_{s(n)}(\X))^2}.\]
Since $\Delta(m,A_{s(n)}(\X))\leq \Delta(m,[0,1]^p)<\infty$, we can use the dominated convergence theorem, and conclude that, using \Cref{P1}, 
\[\lim_{n\to\infty}\E[m_n''(\X)-m(\X)]^2=0.\]
Hence, we have the consistency $m_n'(\X)\xrightarrow{\P}m(\X)$. This means that a (fictive) estimator built upon the empirical partition but where all the values are observed to compute the empirical mean step is consistent. We will use this fact to show $m_n(\X) \xrightarrow{\P} m(\X)$. 

First, consider the case in dimension 1. We use the specific case where a cut $d$ and an assignation $w$ (of the cell $A_{s(n)}$)  leaves all the observed points to the left and assigns all the missing observations to the right. By \Cref{L2}, we know that for $n$ sufficiently large $L_n(A_{s(n)},d,w)\leq\xi$. As already seen in \Cref{eq:CART},
\[L_n(A_{s(n)},d,w)=\frac{\widehat{N}_{obs}(A_{s(n)})\widehat{N}_{miss}(A_{s(n)})}{\widehat{N}(A_{s(n)})\widehat{N}(A_{s(n)})}\parent{\widehat{Y}_{obs}-\widehat{Y}_{miss}}^2.\]
By the convergence of $m_n'(\X)$, we have, in particular that $\widehat{Y}_{obs}\xrightarrow{\P}m(\X)$.

Using \Cref{H2} and the same ideas as in \Cref{L2}, since $p_n^{(h)}\to c^{(h)}<1$, there exists $N\in\N^\star$ such that, with probability at least $1-\rho$ for all $n\ge N$, $\widehat{N}_{obs}(A_{s(n)})/\widehat{N}(A_{s(n)})\ge c$, where $c>0$ is a constant. On the other hand, if $\widehat{N}_{miss}(A_{s(n)})/\widehat{N}(A_{s(n)})\xrightarrow{\P} 0$ then, trivially $m_n(\X)\xrightarrow{\P} m_n'(\X)$ so let us consider the case $\widehat{N}_{miss}(A_{s(n)})/\widehat{N}(A_{s(n)})\ge c'$, hence with probability at least $1-\rho$,
\[\parent{\widehat{Y}_{obs}-\widehat{Y}_{miss}}^2\leq \frac{\xi}{cc'}.\]
This shows that the random variable $\widehat{Y}_{obs}-\widehat{Y}_{miss}$ converges to 0 in probability. Since $\widehat{Y}_{obs}\xrightarrow{\P}m(\X)$, we obtain that $\widehat{Y}_{miss}\xrightarrow{\P}m(\X)$. Finally, using the following formula for $m_n(\X)$
\[m_n(\X)=\frac{\widehat{N}_{obs}(A_{s(n)}(\X))}{\widehat{N}(A_{s(n)}(\X))}\widehat{Y}_{obs}+\frac{\widehat{N}_{miss}(A_{s(n)}(\X))}{\widehat{N}(A_{s(n)}(\X))}\widehat{Y}_{miss},\]
we conclude that $m_n(\X)\xrightarrow{\P}m(\X)$.
For dimension bigger than 1, denote again \(Y^{(j)}=m_j(X^{(j)})+\epsilon^{(j)}\), where \(\epsilon^{(j)}\sim \mathcal{N}(0,\sigma^2/p^2)\) so that we have \(Y\sim\sum_{j=1}^pY^{(j)}\) and define
\[m_n^{(j)}(\X)=\frac{1}{\widehat{N}(A_{s(n)}(\X))}\sum_{i=1}^n Y_i^{(j)}\indicadora{\widehat{\X}_i\in A_{s(n)}(\X)}{}.\]
Consider the cut in the direction $1$ which leaves all the observations where $\M^{(1)}=0$ to the left and assigns the observations where $\M^{(1)}=1$ to the right. We denote $ \widehat{Y}_{obs(1)}$ and $\widehat{Y}_{miss(1)}$ the respective (on the left and on the right) empirical means. By the arguments in dimension 1, we have that $\widehat{Y}_{obs(1)}-\widehat{Y}_{miss(1)} \to 0$ in probability. By definition,
\[\widehat{Y}_{obs(1)}-\widehat{Y}_{miss(1)} = \widehat{Y}_{obs}^{(1)}-\widehat{Y}_{miss}^{(1)}+\left(\sum_{j=2}^p \widehat{Y}_{obs(1)}^{(j)}-\sum_{j=2}^p\widehat{Y}_{miss(1)}^{(j)}\right).\]
Since the random variable $Y^{(j)}$ ($j\neq 1$) is independent of the random variable $\widehat{\X}^{(1)}$ conditionally to $\widehat{\X}\in A_{s(n)}$, the distributions of the two sums on the right hand side are equal. Since each random variable $\widehat{Y}^{(j)}$ converges (see Technical Lemma \ref{TL1}) we conclude that the difference of the two sums converges in probability to 0. Hence, $\widehat{Y}_{obs}^{(1)}-\widehat{Y}_{miss}^{(1)} \to 0$ in probability. This finally shows that $m_n^{(1)}(\X)\xrightarrow{\P}m_1(\X)$. Similarly, we show that for all $j\ge 1$, $m_n^{(j)}(\X)\xrightarrow{\P}m_j(\X)$ and then $m_n(\X)\xrightarrow{\P}m(\X)$ by summation of the $p$ previous convergences which concludes the proof of the Theorem.

\section{Simulation Study}
\label{sec:simulations}

We present a brief simulation study for our proposed method to handle missing values in a regression task. The study compares this approach with other 5 methods to handle missing values with random forests. One of them, corresponding to median imputation before the construction of the random forest is taken as a simple baseline. Three methods correspond to more elaborated imputation algorithms corresponding to the Breiman’s method presented in \citet{Breiman2003}, the improvement suggested by \citet{ishioka2013imputation} (here presented simply as Ishioka’s method) and missForest \citep{stekhoven2011missforest}. Finally, the last method considered in this study corresponds to Missing Incorporated in Attributes (MIA) \citep{twala2008good} which uses the observations with missing values directly in the construction of the recursive trees. Two of these approaches, corresponding to missForest and MIA could be considered as state-of-the-art algorithms to handle missing values through random forests\footnote{Codes to reproduce our results can be found in: \url{https://github.com/IrvingGomez/RandomForestsSimulations}}.

\subsection{Approaches considered in the simulation}
To properly introduce the algorithms considered in this study, we need to define the connectivity between two points in a tree and the proximity matrix of the forest. Let $K_{\Theta,n}(\X,\X')$ be the indicator that $\X$ is in the same final cell that $\X'$ in the tree designed with $\sample$ and the parameter $\Theta$. If $K_{\Theta,n}(\X,\X')=1$ we say that $\X$ and $\X'$ are connected in the tree $m_n(\cdot;\Theta)$. The proximity of two points is the average of times in which they are connected in the forest. Formally, let us define the proximity between $\X$ and $\X'$ in the finite forest, $m_{M,n}(\cdot;\Theta_1,\ldots,\Theta_M)$, as
\[K_{M,n}(\X,\X')=\frac{1}{M}\sum_{k=1}^M K_{\Theta_k,n}(\X,\X')\]

Breiman's and Ishioka's approach operate through imputation of missing values in a recursive way. First, they use the original training data set, $\sample$, to fill the blank spaces in a roughly way. We denote this new data set as $\mathcal{D}_{n,t_1}$. The imputed data set $\mathcal{D}_{n,t_1}$ is used to build a random forest. Then the proximity matrix is used to improve the imputation, resulting in a new data set $\mathcal{D}_{n,t_2}$. The procedure follows iteratively until some stopping rule is applied. Let $K_{M,t_\ell}(i,j)$ be the proximity between $\X_i$ and $\X_j$ at time $t_\ell$. And denote by $\widehat{\X}_{i,t_\ell}$  the imputation of the observation $\X_i$ at time $t_\ell$. We also define $\mathbf{i}_{miss}^{(h)}\subseteq\{1,\ldots,n\}$ as the indexes where $\X^{(h)}$ is missing, and $\mathbf{i}_{obs}^{(h)}=\{1,\ldots,n\}\setminus\mathbf{i}_{miss}^{(h)}$ as the indexes were $\X^{(h)}$ is observed.

\paragraph{Breiman's Approach.}
If $\X^{(h)}$ is a continuous variable, $\widehat{\X}_{j,t_{\ell+1}}^{(h)}$ is the weighted mean of the observed values in $\X^{(h)}$, where the weights are defined by the proximity matrix of the current random forest, that is
\[\widehat{\X}_{j,t_{\ell+1}}^{(h)}=\frac{\sum_{i\in\mathbf{i}_{obs}^{(h)}}K_{M,t_{\ell}}(i,j)\X_i^{(h)}}{\sum_{i\in\mathbf{i}_{obs}^{(h)}}K_{M,t_{\ell}}(i,j)},\quad \begin{array}{l}\ell\geq 1\\
j\in\mathbf{i}_{miss}^{(h)}
\end{array}\]
On the other hand, if $\X^{(h)}$ is a categorical variable, $\widehat{\X}_{j,t_{\ell+1}}^{(h)}$ is given by
\[\widehat{\X}_{j,t_{\ell+1}}^{(h)}=\argmax_{\x\in\mathcal{X}^{(h)}}\sum_{i\in \mathbf{i}_{obs}^{(h)}}K_{M,t_{\ell}}(i,j)\indicadora{\X_i^{(h)}=\x}{},\quad \begin{array}{l}\ell\geq 1\\
j\in\mathbf{i}_{miss}^{(h)}
\end{array}\]
That is, $\widehat{\X}_{j,t_{\ell+1}}^{(h)}$ is the class that maximizes the sum of the proximity considering the observed values in the class.

\paragraph{Ishioka's Approach.}
If $\X^{(h)}$ is a continuous variable, $\widehat{\X}_{j,t_{\ell+1}}^{(h)}$ is the weighted mean of the $k$ nearest neighbors, according to the proximity matrix, over all the values, both imputed and observed. The $k$ closest values are chosen to make more robust the method and avoid values which are outliers.
\[\widehat{\X}_{j,t_{\ell+1}}^{(h)}=\frac{\sum_{\begin{smallmatrix} i\in\text{neigh}_k \\ i\neq j
\end{smallmatrix}}K_{M,t_{\ell}}(i,j)\widehat{\X}_{i,t_{\ell}}^{(h)}}{\sum_{\begin{smallmatrix} i\in\text{neigh}_k \\ i\neq j
\end{smallmatrix}}K_{M,t_{\ell}}(i,j)},\quad \begin{array}{l}\ell\geq 1\\
j\in\mathbf{i}_{miss}^{(h)}
\end{array}\]
For categorical variables, it is not necessary to see only the $k$ closest values because the outliers of $\X$ will have few attention. Meanwhile the proximity with missing values should have more attention, specially when the missing rate is high. Hence, if $\X^{(h)}$ is a categorical variable, $\widehat{\X}_{j,t_{\ell+1}}^{(h)}$ is given by

\[\widehat{\X}_{j,t_{\ell+1}}^{(h)}=\argmax_{\x\in\mathcal{X}^{(h)}}\sum_{i\neq j}K_{M,t_{\ell}}(i,j)\indicadora{\X_{i,t_{\ell}}^{(h)}=\x}{},\quad \begin{array}{l}\ell\geq 1\\
j\in\mathbf{i}_{miss}^{(h)}
\end{array}\]

\paragraph{MissForest.}
This algorithm begins with a rough imputation for the missing values, for example the median of each direction. Then, for each direction $\X^{(h)}$ a random forests is built using all the other directions $\X^{(1)},\ldots,\X^{(h-1)},\X^{(h+1)},\ldots,\X^{(p)}$ and the response, where the random forests is constructed using the observations such that $\M^{(h)}=0$. Then the missing values of $\X^{(h)}$ are predicted with this random forests, these predictions correspond to the imputation of the missing values. These steps are repeated iteratively until a stopping rule is achieved.

\paragraph{Missing Incorporated in Attributes (MIA).}
The Missing Incorporated in Attributes (MIA) approach consists in keeping all the missing values together when a split is performed. Thus, the splits with this approach assign the values according to one of the following rules:

\begin{itemize}
\item $\{\X^{(h)}<z\text{ and }\M^{(h)}=1\}$ versus $\{\X^{(h)}\geq z\}$.
\item $\{\X^{(h)}<z\}$ versus $\{\X^{(h)}\geq z\text{ and }\M^{(h)}=1\}$.
\item $\{\M^{(h)}=0\}$ versus $\{\M^{(h)}=1\}$.
\end{itemize}

\subsection{Parameters considered and results}
The regression function is the so-called “friedman1” \citep{friedman1991multivariate}, which has been used in previous simulation studies \citep{breiman1996bagging,rieger2010random,josse2019consistency,friedberg2020local}, given by \[m(\x)=10\sin\parent{\pi\x^{(1)}\x^{(2)}}+20\parent{\x^{(3)}-0.5}^2+10\x^{(4)}+5\x^{(5)}.\]
Following the schema presented by \citet{rieger2010random}, we simulate $\X$ as a uniformly distributed variable on $[0,1]^5$ and introduce missing values in $\X^{(1)},\X^{(3)}$ and $\X^{(4)}$, considering an MCAR data-missing mechanism. We create one testing data set and 100 training data sets. Each training data set contains 200 observations and the testing data set contains 2000 observations. This amount of data is to have an appropriate approximation to the MSE \[\E_{\X|\sample}\corchet{m_{M,n}(\X)-m(\X)}^2.\]

For each training data set $\X^{(1)}$ presents 20\% of missing data, $\X^{(3)}$ has 10\% of missing data, while for  $\X^{(4)}$ we have decided to change the fraction of missingness varying between 5\%, 10\%, 20\%, 40\%, 60\%, 80\%, 90\% and 95\%.  We do not introduce missing values in the testing data set. A random forest is built for each one of the 100 training data sets with missing values and without missing values (which are used as benchmark), using $M=50$ trees, that has been seen by simulation to be sufficient to stabilize the error in the case of complete training data sets. For the remain parameters we use the default values in the regression mode of the \texttt{R} package \texttt{randomForests}, the parameter $\mathtt{mtry}$ is set to $\lfloor p/3\rfloor$, we have sampled without replacement, so $a_n$ is set to $\lceil 0.632 n\rceil$ and $\texttt{nodesize}$ is set to 5.

\begin{figure}[H]
\begin{center}
\includegraphics[width=\textwidth]{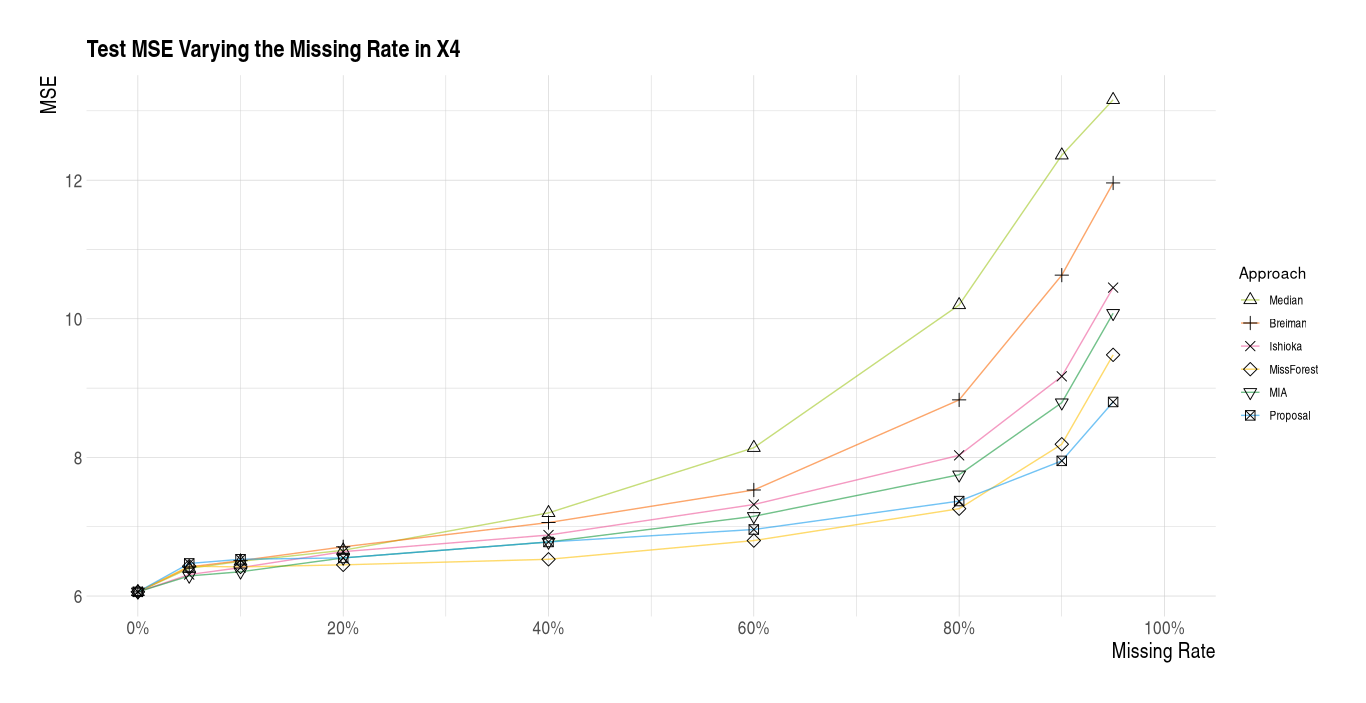}
\caption{Average MSE for the testing data set for each percentage of missingness}
\label{fig:simulations}
\end{center}
\end{figure}

\Cref{fig:simulations} presents the average MSE over the training data sets for each percentage of missing data considered. Median imputation could be taken as a bound of the minimum expected behavior for a method that attempts to estimate the regression function with missing values. We observe that missForest and the approach here presented generate estimators with the lowest MSE regardless of the percentage of missingness. Important differences between the performance of the methods become clear with the increasing percentage of missing values, especially when this value is over 60\%. Moreover, we can see the advantage of searching for the best assignation when there is a high percentage of missing values with our proposal outperforming all the other algorithms. While all the methods present a similar MSE when this percentage is smaller than 20\%.

Another interesting observation corresponds to the relation between our procedure and Missing Incorporated in Attributes (MIA), which assigns all the observations where the split variable is missing to the same node. Our method goes further than MIA, since it does not only consider these assignations, but the best possible partial imputation for the observations with missing values. Furthermore, note that the algorithm proposed in this work has the advantage of not having to calculate extra structures, like the proximity matrix. The work presented in \citet{gomez2021rfsimulation} is devoted to the computational study of the algorithm here proposed, it presents results of an extensive simulation study considering different data-missing mechanisms beyond MCAR.

\printbibliography

\end{document}